\documentclass[11pt]{article}
\usepackage{amsfonts}
\usepackage{graphics}
\usepackage{indentfirst}
\usepackage{color}
\usepackage{cite}
\usepackage{tikz}
\usepackage{latexsym}
\usepackage[paper=a4paper, left=2.35cm, right=2.35cm, top=1.8cm, bottom=1.8cm, headheight=5.5pt, footskip=0.8cm, footnotesep=0.8cm, centering, includefoot]{geometry}
\usepackage{amsmath}
\allowdisplaybreaks
\usepackage{amssymb}
\usepackage[colorlinks, linkcolor=red]{hyperref}
\hypersetup{urlcolor=red, citecolor=red}
\usepackage[dvips]{epsfig}
\usepackage{amscd}

\usepackage{amsthm}
\usepackage{mathrsfs}
\usepackage{verbatim}

\newtheorem{theorem}{Theorem}[section]
\newtheorem{lemma}[theorem]{Lemma}
\newtheorem{coro}[theorem]{Corollary}

\newtheorem{conjecture}[theorem]{Conjecture}
\newtheorem{prop}[theorem]{Proposition}
\theoremstyle{definition}
\newtheorem{defn}[theorem]{Definition}

\newtheorem{exmp}[theorem]{Example}

\newcommand{\delete}[1]{}

\makeatletter
\@addtoreset{equation}{section}
\makeatother
\makeatletter
\@addtoreset{equation}{section}
\makeatother

\title{Separable elements and splittings in Weyl groups of type $B$
\thanks{This work was partially supported by the National Natural Science Foundation of China (Grant No. 12471020, 12071377),
the Natural Science Foundation of Chongqing (Grant No. CSTB2023NSCQ-MSX0706), and the Fundamental Research Funds for the Central Universities (Grant No. SWU-XDJH202305).}}

\author{Ming Liu,\ Houyi Yu {\thanks{Corresponding author. \newline\hspace*{1.6em} E-mail addresses:
liuming2647@163.com (M. Liu),  yuhouyi@swu.edu.cn (H. Yu).}}\date{}\\
\footnotesize School of Mathematics and Statistics,
Southwest University, Chongqing 400715, P. R. China}
\begin{document}
\maketitle

\begin{abstract}
Separable elements in Weyl groups are generalizations of the well-known class of separable permutations in symmetric groups.
Gaetz and Gao showed that for any pair $(X,Y)$ of subsets of the symmetric group $\mathfrak{S}_n$,
the multiplication map $X\times Y\rightarrow \mathfrak{S}_n$ is a splitting (i.e., a length-additive bijection) of $\mathfrak{S}_n$ if and only if $X$ is the generalized quotient of $Y$ and
$Y$ is a principal lower order ideal in the right weak order generated by a separable element.
They conjectured this result can be extended to all finite Weyl groups. In this paper, we classify all separable
and minimal non-separable signed permutations in terms of forbidden patterns and confirm the conjecture of Gaetz and Gao for Weyl groups of type $B$.
\end{abstract}


%2020 \textit{Mathematics Subject Classification}. 76A15; 76N10; 35Q35.

%\tableofcontents

\section{Introduction}\label{sec:int}

Let $\mathfrak{S}_n$ denote the symmetric group of all permutations of the set $[n]=\{1,2,\ldots,n\}$.
A permutation $w=w_1w_2\cdots w_n$ in $\mathfrak{S}_n$ is said to be separable if it avoids the patterns $3142$ and $2413$,
that is, there do not exist indices $i_1<i_2<i_3<i_4$ such that $w_{i_2}<w_{i_4}<w_{i_1}<w_{i_3}$ or $w_{i_3}<w_{i_1}<w_{i_4}<w_{i_2}$.
Separable permutations were formally introduced by Bose, Buss, and Lubiw \cite{BBL98} to deal with the pattern matching problem for permutations,  although these permutations can be tracked back to the much earlier work of Avis and Newborn \cite{AN81} on pop-stacks.
An account of the research on these permutations can be found in  \cite{Kit11,Vat15}.
Separable permutations possess many amazing properties. For instance, they
can be completely constructed recursively and represented as signed Schr\"oder trees \cite{BBFGP18}, are enumerated by the large
Schr\"oder numbers \cite{Kre00,Wes95}, and have interesting connections with permutation pattern matching \cite{BBL98,NRV16} and bootstrap percolation \cite{Sha91}.

In 2020, Gaetz and Gao \cite{GG20aam} generalized separable permutations in symmetric groups to separable elements in arbitrary finite Weyl groups.
These elements play an essential role in studying Weyl groups.
Gaetz and Gao \cite{GG20am, GG20aam} showed that separable elements have a characterization by root system pattern avoidance in the sense of Billey and Postnikov \cite{BP05} and can induce certain algebraic decompositions of the Weyl groups.
More recently, Gossow and Yacobi \cite{GY23} found that separable elements of simply-laced Weyl groups act on the canonical bases of the representations of the corresponding groups by bijection up to lower-order terms.

Separable elements are closely related to the algebraic structure of Weyl groups.
Generalized quotients were introduced by Bj\"{o}rner and Wachs \cite{BW88} for Coxeter groups.
Recall that given a subset $U$ of a Coxeter group $W$, the subset
\begin{align*}
   W/U=\{w\in W \mid \ell(wu)=\ell(w)+\ell(u)\ \text{for all}\ u\in U\}
\end{align*}
is called a \emph{generalized quotient}.
The name is motivated by the fact that if $U$ is a parabolic subgroup, say $W_J$, then its generalized quotient $W/W_J$ is exactly the parabolic quotient $W^J$.
Bj\"{o}rner and Wachs \cite{BW88} showed that for a finite Coxeter group $W$, the generalized quotient $W/U$ is always an interval $[e, w_0u^{-1}]_L$ in the left weak order, where $u$
is the least upper bound of $U$ in the right weak order, and $w_0$ is the longest element of $W$.

A pair $(X,Y)$ of subsets of a Weyl group $W$ is called a \emph{splitting} of $W$ if the multiplication map $X\times Y\rightarrow W$ is length-additive (i.e., $\ell(xy)=\ell(x)+\ell(y)$ for all $x\in X$ and $y\in Y$) and bijective. For instance, $(W^J,W_J)$ is a splitting of $W$.

For each separable element of a finite Weyl group $W$, Gaetz and Gao \cite[Theorem 3]{GG20am} gave a splitting  to $W$ by using generalized quotients. More precisely,
given any element $u$ in a Weyl group $W$, if $U$ is the interval $[e,u]_R$ in the right weak order with $u$ separable, then $(W/U,U)$ is a splitting of $W$.
Moreover, they showed for the symmetric group that any splitting of this form can be induced by a separable permutation, and conjectured this extends to arbitrary Weyl groups:
\begin{conjecture}[\cite{GG20am}, Conjecture 2]\label{conj2}
Let $W$ be a finite Weyl group, and let $U=[e,u]_R\subseteq W$. Then $(W/U,U)$ is a splitting of $W$ if and only if $u$ is separable.
\end{conjecture}
The main result of this paper is the following theorem, which confirms Conjecture \ref{conj2} for Weyl groups of type $B$.
\begin{theorem}\label{maintheorem}
Let $W$ be a Weyl group of type $B$, and let $U=[e,u]_R\subseteq W$. Then $(W/U,U)$ is a splitting of $W$ if and only if $u$ is separable.
\end{theorem}

This paper is organized as follows. Section \ref{sec:Pre} introduces necessary preliminary results on Weyl groups and separable elements.
In Section \ref{Sec:separableelementsinBn}, we present the standard combinatorial descriptions of Weyl groups of type $B$,
and then give alternative characterizations of separable elements and minimal non-separable elements in this group in terms of forbidden patterns. The proof of Theorem \ref{maintheorem} is concluded in Section \ref{sec:Proofofthemainythm}.

\section{Preliminaries}\label{sec:Pre}

In this section, we fix some notation and collect basic facts concerning posets, Weyl groups, and separable elements, and referring the reader to \cite{BB05,Hum90,Sta12} for further details.

Given any integers $m$ and $n$, let $[m,n]$ denote the set $\{k\in\mathbb{Z}\mid m\leq k\leq n\}$ whenever $m\leq n$ and the empty set $\emptyset$ { otherwise}. We simply write
$[n]=\{1,2,\ldots,n\}$ and $[n]_{\pm}=\{1,2,\ldots,n\}\cup\{\underline{1},\underline{2},\ldots,\underline{n}\}$, where underlined numbers are negative.
Unless otherwise stated, $n$ will denote a nonnegative integer throughout the paper.
For any number $x$, let $|x|$ denote the absolute value of $x$.
We use $\uplus$ to denote the disjoint union and $\# X$ to denote
the cardinality of a set $X$.

\subsection{Rank-generating functions of graded posets}

A sequence $a_0, a_1, \ldots,a_n$ of real numbers is \emph{symmetric} if $a_i=a_{n-i}$ for all $i\in[0,n]$, and is
\emph{unimodal} if
$$a_0\leq\cdots\leq a_{m-1}\leq a_m\geq a_{m+1}\geq\cdots\geq a_n$$
for some $m\in[0,n]$.
Similarly, a polynomial $f(q)=\sum_{i=0}^n a_iq^i$ is said to be \emph{symmetric} (sometimes called \emph{{ palindromic}}), if the sequence $a_0, a_1, \ldots,a_n\neq0$ of coefficients is symmetric, and is \emph{unimodal} if this sequence is unimodal.
Given a polynomial $f(q)$, we write $[q^i]f$ to represent the coefficient of $q^i$ in $f$.
It is well known that the product of two symmetric and unimodal polynomials with nonnegative coefficients is still symmetric and unimodal, see, for example, \cite[Proposition 1]{Sta86}.

Let $P$ be a poset. For each element $x$ of $P$, let $V_x=\{y\in P\mid y\geq x\}$ denote the principal upper order ideal generated by $x$. Similarly, $\Lambda_x=\{y\in P\mid y\leq x\}$ denotes the principal lower order ideal generated by $x$.
If $P$ is a graded poset of rank $n$ and has $p_i$ elements of rank $i$, then
the polynomial $P(q)=\sum_{i=0}^np_iq^i$ is called the \emph{rank-generating function} of $P$.
We say a graded poset $P$ is \emph{rank-symmetric} if $P(q)$ is symmetric, and is \emph{rank-unimodal} if $P(q)$ is unimodal.

\subsection{Root systems and Weyl groups}\label{subsec:RSinWG}
Let $\Phi$ be a crystallographic root system with a simple system $\Delta$ and the corresponding positive system $\Phi^+$, and let $W=W(\Phi)$ be the Weyl group generated by \emph{simple reflections} $S=\{s_\alpha\mid \alpha\in\Delta\}$.
Each element $w$ in $W$ can be written as a product of simple reflections $w=s_1s_2\cdots s_r$. If $r$ is minimal, then such an expression is called \emph{reduced}, and $r$ is called the \emph{length} of $w$, denoted $\ell(w)=r$. The \emph{inversion set} of $w$ is $I_{\Phi}(w)=\{\alpha\in\Phi^+\mid w\alpha\in \Phi^-\}$.
It's well known that $\ell(w)=\#I_{\Phi}(w)$. There exists a unique element $w_0$ that satisfies $I_{\Phi}(w_0)=\Phi^+$, called the \emph{longest element}.
For convenience, we define a partial order on $\Phi^{+}$ such that $\alpha\leq \beta$ if $\beta-\alpha$ is a nonnegative (integral) linear combination of the simple roots.

For $J\subseteq \Delta$, we denote by $W_J$ the \emph{parabolic subgroup} of $W$ generated by the reflections associated to the roots in $J$. The \emph{parabolic quotient} of $W_J$ in $W$ is
$$W^J=\{w\in W\mid \ell(w)<\ell(ws_{\alpha})\ {\rm for\ all}\ \alpha\in J\},$$
which intersects each left coset of $W_J$ in its unique element of minimum length. Then every $w$ in $W$ has a unique factorization $w=w^J\cdot w_J$ such that $w^J\in W^J$ and $w_J\in W_J$. For this factorization, we have $\ell(w)=\ell(w^J)+\ell(w_J)$.

The \emph{left weak order} on $W$ is defined by $u\leq_Lw$ if $\ell(w)=\ell(u)+\ell(wu^{-1})$, or equivalently $I_{\Phi}(u)\subseteq I_{\Phi}(w)$.
The \emph{right weak order} on $W$ is defined analogously by $u\leq_Rw$ if  $\ell(w)=\ell(u)+\ell(u^{-1}w)$.
The two weak orders are isomorphic via the map $w\mapsto w^{-1}$.
Each weak order makes the group $W$ a graded poset ranked by the length function, and $W$ is indeed a lattice, with minimal element $e$ and maximal element $w_0$.
So we have $\Lambda^L_w=[e,w]_L$ and  $V^L_w=[w,w_0]_L$ for any $w\in W$.
If $W$ has rank $n$, then
\begin{align}\label{eq:W(q)}
  W(q)=\sum_{w\in W}q^{\ell(w)}=\prod_{i=1}^n\frac{q^{d_i}-1}{q-1},
\end{align}
where the $d_i$ are integer invariants called the \emph{degrees} of $W$.

A root system $\Phi$ is \emph{irreducible} if it cannot be properly partitioned into $\Phi_1\uplus\Phi_2$
such that $s_{\alpha}s_{\beta}=s_{\beta}s_{\alpha}$ for all $\alpha\in\Phi_1$ and $\beta\in\Phi_2$.
It is well known that (see, for example, \cite{Hum90}) irreducible crystallographic root systems can be
classified into $4$ infinite families of types $A_n$, $B_n$, $C_n$, $D_n$ and exceptional types
$E_6$, $E_7$, $E_8$, $F_4$, and $G_2$.

Let $\Phi$ be a crystallographic root system.
A subset $\Phi'$ of $\Phi$ is called a \emph{subsystem} if $\Phi'=\Phi\cap U$, where $U$ is a linear subspace of ${\rm span}(\Phi)$. It is easy to see that $\Phi'$ is a root system itself
in the subspace spanned by $\Phi'$.
The notation of pattern avoidance in Weyl groups was first introduced by Billey and Postnikov \cite{BP05}.
For $w\in W(\Phi)$, we say $w$ \emph{contains the pattern} $(w', \Phi')$ if $I_{\Phi}(w)\cap U=I_{\Phi'}(w')$; we write $w|_{\Phi'}=w'$ in this case. If $\Phi'$ is spanned by $J\subset\Delta$, then $W(\Phi')$ can be identified  with $W_J$ naturally, so $w|_{\Phi'}=w_J$. We say $w$ \emph{avoids} $(w',\Phi')$ if it does not contain any pattern isomorphic to $(w',\Phi')$.

\subsection{Separable elements in Weyl groups}
Let $\Phi$, $\Delta$, and $W(\Phi)$ be as in Subsection \ref{subsec:RSinWG}.
The following notation of separable elements in general Weyl groups was first defined in a recursive manner by Gaetz and Gao \cite{GG20am,GG20aam}, generalizing the notion of separable permutations \cite{AN81}.

\begin{defn}\label{separabledef}
Let $w\in W(\Phi)$. Then $w$ is \emph{separable} if one of the following holds:
\begin{itemize}
  \item[(S1)] $\Phi$  is of type $A_1$;
  \item[(S2)] $\Phi=\bigoplus\Phi_i$ is reducible and $w|_{\Phi_i}$ is separable for each $i$;
  \item[(S3)] $\Phi$ is irreducible and there exists a \emph{pivot} $\alpha_i\in \Delta$ such that $w|_{\Phi_J}\in W(\Phi_J)$ is separable, where $\Phi_J$ is generated by $J=\Delta\backslash\{\alpha_i\}$ and such that either
\begin{align*}
   \{\beta\in \Phi^+\mid\beta\geq\alpha_i\}\subset I_{\Phi}(w),\quad \text{or}\quad
   \{\beta\in \Phi^+\mid\beta\geq\alpha_i\}\cap I_{\Phi}(w)=\emptyset.
\end{align*}
%We say $\alpha_i$ is a \emph{full pivot} in the first case, and an \emph{empty pivot} in the second.
\end{itemize}
\end{defn}
Note that the separability depends not only on $w$ but also on $\Phi$, so we sometimes write the separable element $w$ as $(w,\Phi)$ for clarity.
In \cite{GG20aam}, Gaetz and Gao gave an equivalent characterization for separable elements in terms of root system pattern avoidance.

\begin{lemma}[\cite{GG20aam}, Theorem 5.3]\label{patternavoid}
An element $w\in W(\Phi)$ is separable if and only if $w$ avoids the following root system patterns:
\begin{enumerate}
  \item\label{patternavoidA3} the patterns corresponding to the permutations $3142$ and $2413$ in the Weyl group of type $A_3$,
  \item\label{patternavoidB2} the two patterns of length two in the Weyl group of type $B_2$, and
  \item\label{patternavoidG2} the six patterns of length two, three and four in the Weyl group of type $G_2$.
\end{enumerate}
\end{lemma}

The following remarkable property of separable elements was first proved by Wei \cite{We12} for symmetric groups, and later by Gaetz and Gao \cite{GG20aam}
for arbitrary Weyl groups.

\begin{lemma}[\cite{GG20aam}, Theorem 3.9]\label{rank-symmetric}
Let $W$ be a Weyl group, and let $w\in W$ be separable. Then the left weak order order ideals $V_w^L$ and $\Lambda_w^L$ are both rank-symmetric and rank-unimodal, and
$\Lambda^L_w(q) V_w^L(q)=W(q)$.
\end{lemma}

\section{Separable elements in Weyl groups of type $B$}\label{Sec:separableelementsinBn}
In preparation for our main results,
we present in this section the standard combinatorial descriptions of the Weyl group of type $B_n$ ($n\geq2$).
Then we classify all separable elements and minimal non-separable elements in this group in terms of signed permutation pattern avoidance.

We concern mainly with the Weyl group of type $B_n$, so throughout the remainder of this paper, we adopt the following conventions for the root system of type $B_n$:
$$\Phi=\{\pm e_i\mid i\in[n]\}\cup \{\pm e_i\pm e_j\mid {1\leq i<j\leq n}\}\subset\mathbb{R}^n$$
with positive roots
$$\Phi^+=\{\pm e_i+ e_j \mid 1\leq i<j \leq n\}\cup\{e_i\mid i\in[n]\}$$
and simple roots $$\Delta=\{\alpha_0=e_1\}\cup\{\alpha_i=-e_i+e_{i+1}\mid i\in[n-1]\}.$$
Here the $e_i$ denotes the $i$-th standard basis vector in $\mathbb{R}^n$.
Note that we have chosen a different simple system from \cite{GG20am,GG20aam}.
\delete{Figure $1$ illuistrates the Coxeter graph of $B_n$.
\begin{center}
\begin{tikzpicture}\label{Figure1}
\node (alpha0) at (-3,0) {{$\bullet$}};
\node (alpha1) at (-1.5,0) {{$\bullet$}};
\node (alpha2) at (0,0) {{$\bullet$}};
\node (alpha) at (0.8,0) {{{$\cdots$}}};
\node (alpha{n-2}) at (1.5,0) {{$\bullet$}};
\node (alpha{n-1}) at (3,0) {{$\bullet$}};
\node (alpha0) at (-2.25,0.3) {{$4$}};
\node (alpha0) at (-3,-0.4) {{$\alpha_0$}};
\node (alpha1) at (-1.5,-0.4) {{$\alpha_1$}};
\node (alpha2) at (0,-0.4) {{$\alpha_2$}};
\node (alpha{n-2}) at (1.5,-0.4) {{$\alpha_{n-2}$}};
\node (alpha{n-1}) at (3,-0.4) {{$\alpha_{n-1}$}};
\draw[black, thick] (-3,0)--(0,0);
\draw[solid,line width=0.8pt] (1.5,0)--(3,0);
\node(Figure 1.)   at (0,-1.2) {Figure $1.$ The Coxeter graph of $B_n$.};
\end{tikzpicture}
\end{center}
}

Following \cite[Chapter 8]{BB05}, with the above choice of simple system, the Weyl group of type $B_n$  ($n\geq2$) has a combinatorial interpretation as the group $\mathfrak{B}_n$ of signed permutations.
Recall that a \emph{signed permutation} $w$ of $[n]$ is a bijection of the set $[n]_{\pm}$ onto itself such that $w(\underline i)=\underline{w(i)}$ for all $i\in[n]$.
%The set of all signed permutations of $[n]$ naturally forms the $n$-th \emph{hyperoctahedral group} which denoted in the literature by.
For an element $w\in \mathfrak{B}_n$, we use both the one-line notation $w=w_1w_2\cdots w_n$ and the disjoint cycle notation, where $w_i=w(i)$ for $i\in[n]$.
As a set of generators for $\mathfrak{B}_n$, we take $S^B=\{s_0,s_1,\ldots,s_{n-1}\}$, where $s_0=(1,\underline 1)$ and $s_{i}=(i,i+1)$ for all $i\in[n-1]$.
The Weyl group of type $A_{n-1}$ (i.e., the symmetric group $\mathfrak{S}_n$) is the subgroup in which all elements have no negative entries.

\begin{exmp}
According to Lemma \ref{patternavoid}, the separable elements in $\mathfrak{B}_2$ are
$12$, $\underline{1}2$, $21$, $1\underline{2}$, $\underline{2}\,\underline{1}$, and $\underline{1}\,\underline{2}$,
whose inversion sets are
\begin{align*}
\emptyset,\quad \{\alpha_0\},\quad \{\alpha_1\}, \quad\{\alpha_1,\alpha_0+\alpha_1,2 \alpha_0+\alpha_1\}, \quad\{\alpha_0,\alpha_0+\alpha_1,2 \alpha_0+\alpha_1\}, \quad\text{and}\quad \Phi^{+},
\end{align*}
respectively.
\end{exmp}

Let $w=w_1w_2\cdots w_n$ be an element of $\mathfrak{B}_n$. A positive integer $i\in [n]$ is a \emph{negative index} of $w$ if $w_i<0$. A pair $(i,j)\in[n]^2$ with $i<j$ is called an \emph{inversion} of $w$ if $w_i>w_j$,
and called a \emph{negative sum pair} of $w$ if $w_i+w_j<0$.
Define
\begin{align*}
  {\rm Neg}(w)&=\{i\in[n]\mid w_i<0\},\\
  {\rm Inv}(w)&=\{(i,j)\in [n]^2\mid i<j, w_i>w_j \},\\
  {\rm Nsp}(w)&=\{(i,j)\in [n]^2\mid i<j, w_i+w_j<0\}.
\end{align*}
If we write $-e_i=e_{\underline{i}}$ for $i\in[n]$, then any element $w\in \mathfrak{B}_n$ acts on $\mathbb{R}^n$ by $we_{i}=e_{w_i}$.
Hence, ${\rm Neg}(w)=\{i\in[n]\mid e_i\in I_{\Phi}(w)\}$,
${\rm Inv}(w)=\{(i,j)\in [n]^2\mid -e_i+e_j\in I_{\Phi}(w)\}$,
and ${\rm Nsp}(w)=\{(i,j)\in [n]^2\mid i<j, e_i+e_j\in I_{\Phi}(w)\}$.

Since $u\leq_Lv$ if and only if $I_{\Phi}(u)\subseteq I_{\Phi}(v)$, we get the following alternative characterization of the left weak order on $\mathfrak{B}_n$,
 which is useful for our later discussion.
%Since in any Weyl groups, $u\leqslant_Lw$ if and only if $I_{\Phi}(u)\subseteq I_{\Phi}(w)$, see \cite[Proposition 2.1]{HL16}, we have the following result.
\begin{lemma}\label{leftweakorder}
Let $u,v\in \mathfrak{B}_n$. Then $u\leq_Lv$ if and only if ${\rm Inv}(u)\subseteq {\rm Inv}(v)$, ${\rm Neg}(u)\subseteq {\rm Neg}(v)$, and ${\rm Nsp}(u)\subseteq {\rm Nsp}(v)$.
\end{lemma}

The length function $\ell$ can be { combinatorially} described as
\begin{align}\label{eq:eqell(w)}
  \ell(w)=\# {\rm Neg}(w)+ \# {\rm Inv}(w)+\#  {\rm Nsp}(w)=\# {\rm Inv}(w)-\sum_{\{i\in[n]\mid w_i<0\}}w_i.
\end{align}

Let $a=a_1a_2\cdots a_n$ be a sequence of $n$ nonzero real numbers. The \emph{standard permutation} st$(a)$ of $a$ is the unique permutation $w\in \mathfrak{S}_n$ defined by
\begin{align*}
  w_i<w_j~~~\Leftrightarrow~~~ a_i\leq a_j
\end{align*}
for all $i$, $j$ with $1\leq i<j\leq n$, while the \emph{standard signed permutation} sts($a$) of $a$ is the unique signed permutation
$w\in \mathfrak{B}_n$ such that ${{\rm Neg}}(w)={{\rm Neg}}(a)$ and
\begin{align*}
  |w_i|<|w_j|~~~\Leftrightarrow~~~ |a_i|\leq |a_j|
\end{align*}
for all $i$, $j$ with $1\leq i<j\leq n$. Here ${\rm Neg}(a)=\{i\in[n]\mid a_i<0\}$.
\begin{exmp}\label{exmp31422413}
It is straightforward to check that
\begin{align*}
\{w\in \mathfrak{B}_4\mid {\rm st}(w)=3142\}
=\,&\{3142,\underline{2}\,\underline{4}\,\underline{1}\,\underline{3},3\underline{1}42,3\underline{2}41,2\underline{3}41,2\underline{4}
31,3\underline{2}4\underline{1},2\underline{3}4\underline{1},\\
&\hspace{20mm}2\underline{4}3\underline{1},1\underline{3}4
\underline{2},1\underline{4}3\underline{2},1\underline{4}2\underline{3},\underline{2}\,\underline{4}1\underline{3},\underline{1}\,
\underline{4}2\underline{3},\underline{1}\,\underline{4}3\underline{2},\underline{1}\,\underline{3}4\underline{2}\}
\end{align*}
and
\begin{align*}
\{w\in \mathfrak{B}_4\mid {\rm st}(w)=2413\}
=\,&\{2413,\underline{3}\,\underline{1}\,\underline{4}\,\underline{2},13\underline{4}2,14\underline{3}2,14\underline{2}3
,24\underline{1}3,\underline{3}2\underline{4}1,\underline{2}3\underline{4}1,\\
&\hspace{20mm}\underline{2}4\underline{3}1,\underline{1}3
\underline{4}2,\underline{1}4\underline{3}2,\underline{1}4\underline{2}3,\underline{3}1\underline{4}\,\underline{2},\underline{3}2
\underline{4}\,\underline{1},\underline{2}3\underline{4}\,\underline{1},\underline{2}4\underline{3}\,\underline{1}\}.
\end{align*}
\end{exmp}

Let $J=\Delta\backslash\{\alpha_{p_1}, \alpha_{p_2}, \ldots, \alpha_{p_k}\}$, where $0\leq p_1< p_2< \cdots<p_k\leq n-1$. Then the parabolic subgroup generated by $J$, denoted $\mathfrak{B}_J$, has the form
\begin{align*}
\mathfrak{B}_J=\mathfrak{B}_{p_1}\times\mathfrak{S}_{p_2-p_1}\times\cdots\times\mathfrak{S}_{n-p_k} .
\end{align*}
The set of minimal left coset representatives of $\mathfrak{B}_J$ is
\begin{align*}
\mathfrak{B}^J=\{w\in \mathfrak{B}_n\mid 0<w_1<\cdots <w_{p_1}, w_{p_1+1}<\cdots <w_{p_2}, \ldots, w_{p_{k}+1}<\cdots<w_{n}\}.
\end{align*}
Then for any $w\in \mathfrak{B}_n$, we have
\begin{align*}
w_J={\rm sts}(w_1\cdots w_{p_1})\times {\rm st}(w_{p_1+1}\cdots w_{p_2})\times\cdots\times {\rm st}(w_{p_{k}+1}\cdots w_{n}).
\end{align*}
Moreover, $w^J$ can be obtained from $w$ by first rearranging the elements of $\{|w_1|, \ldots, |w_{p_1}|\}$ in increasing order in the places $1,\ldots,p_1$, and then rearranging the elements of $\{w_{p_i+1}, \ldots, w_{p_{i+1}}\}$ in increasing order in the places $p_i+1,\ldots,p_{i+1}$ for $i\in[k]$, where $p_{k+1}=n$.

The notion of signed permutation pattern avoidance can be defined in terms of the function ${\rm sts}$.
\begin{defn}
Let $w\in \mathfrak{B}_n$ and $u\in \mathfrak{B}_m$ be signed (possibly unsigned) permutations, where $m$ and $n$ are positive integers with $m\leq n$. We say $w$ \emph{contains the pattern} $u$ if in the one-line notation, there exist positive integers
$1\leq i_1<i_2<\cdots<i_m\leq n$ such that ${\rm sts}(w_{i_1}w_{i_2}\cdots w_{i_m})=u$. Otherwise, we say $w$ \emph{avoids the pattern} $u$.
\end{defn}
%Note that if $u$ has no negative entries, then $u$ can be also viewed as a type $A$ pattern. In this case,  We say $w$ \emph{contains the type $A$ pattern} $u$ if there exist positive %integers $1\leq i_1<i_2<\cdots<i_m\leq n$ such that ${\rm st}(w_{i_1}w_{i_2}\cdots w_{i_m})=u$.
%Nearly every pattern treated in what follows is of type $B$, and we shall generally omit the adjective.

Applying Lemma \ref{patternavoid} to $\mathfrak{B}_n$, we obtain a characterization of separable elements of $\mathfrak{B}_n$ by signed permutation pattern avoidance.

\begin{lemma}\label{sepBnavoidpattern}
Let $w=w_1w_2\cdots w_n$ be an element of $\mathfrak{B}_n$. Then $w$ is separable if and only if it avoids the patterns $\underline{2}1$, $2\underline{1}$, $3142$, $2413$, $\underline{3}\,\underline{1}\,\underline{4}\,\underline{2}$, and $\underline{2}\,\underline{4}\,\underline{1}\,\underline{3}$.
\end{lemma}
\begin{proof}It suffices to show that $w$ is non-separable if and only if $w$ contains at least one of the patterns $\underline{2}1$, $2\underline{1}$, $3142$, $2413$, $\underline{3}\,\underline{1}\,\underline{4}\,\underline{2}$, and $\underline{2}\,\underline{4}\,\underline{1}\,\underline{3}$.

If $w$ contains the patterns $\underline{2}1$ or $2\underline{1}$, then there exist $i<j$ in $[n]$ such that ${\rm sts}(w_iw_j)=\underline{2}1$ or $2\underline{1}$.
Let $U$ be the subspace of ${\rm span}(\Phi)$ generated by $\{e_i,-e_i+e_j\}$. Then
$$\Phi\cap U=\{\pm e_i,\pm e_j,{\pm e_i\pm e_j} \},$$
which is isomorphic to the root system of type $B_2$. Clearly, $I_\Phi(w)\cap U=\{e_i,\, e_i+e_j\}$ or $\{e_j,\, -e_i+e_j\}$, from which we see that
$\ell(w|_{\Phi\cap U})=2$.
Hence $w| _{\Phi\cap U}$ is a pattern of type $B_2$ with length two appearing in $w$, and in view of Lemma \ref{patternavoid}\eqref{patternavoidB2} $w$ is non-separable.

If $w$ contains one of the patterns $3142$, $2413$, $\underline{3}\,\underline{1}\,\underline{4}\,\underline{2}$, and $\underline{2}\,\underline{4}\,\underline{1}\,\underline{3}$, then there exist $i<j<k<l$ in $[n]$ such that ${\rm sts}(w_iw_jw_kw_l)\in\{3142, 2413, \underline{3}\,\underline{1}\,\underline{4}\,\underline{2}, \underline{2}\,\underline{4}\,\underline{1}\,\underline{3}\}$.
Let $V={\rm span}\{-e_i+e_j,\, -e_j+e_k,\, -e_k+e_l\}$.
Then $\Phi\cap V$ is  isomorphic to the root system of type $A_3$, and $I_\Phi(w)\cap V=\{-e_i+e_j,\, -e_i+e_l,\, -e_k+e_l\}$ or $\{-e_i+e_k,\, -e_j+e_k,\, -e_j+e_l\}$.
So $w$ contains the pattern $w|_{\Phi\cap V}$, which is a pattern corresponding to the permutations $3142$ and $2413$ in the Weyl group of type $A_3$.
According to Lemma \ref{patternavoid}\eqref{patternavoidA3}, $w$ is non-separable.

Conversely, suppose $w$ is non-separable. Since the ratio of the lengths of the long and short roots are $\sqrt{2}$ in type $B_n$ and $\sqrt{3}$ in type $G_2$, the patterns from Lemma \ref{patternavoid}\eqref{patternavoidG2} never appear in any signed permutation. It follows from  Lemma \ref{patternavoid} that at least one of the following conditions holds:
\begin{enumerate}
\item[(i)]\label{item:ibn}  there exists a subspace $V\subseteq {\rm span}(\Phi)$ such that $\Phi\cap V$ is isomorphic to the root system of type $A_3$ with simple roots $\alpha,\, \beta,\, \gamma \in \Phi^+$ (in this order), and $I_\Phi(w) \cap V=\{\alpha,\, \alpha+\beta+\gamma,\, \gamma\}$ or $\{\alpha+\beta,\, \beta,\, \beta+\gamma\}$;
\item[(ii)]\label{item:iibn2} there exists a subspace $U\subseteq {\rm span}(\Phi)$ such that $\Phi\cap U$ is isomorphic to the root system of type $B_2$, and $\ell(w| _{\Phi\cap U})=2$.		 
\end{enumerate}

For the case (ii), there exist $\alpha,\, \beta$ in $\Phi^+$ as simple roots such that the subsystem ${\rm span}\{\alpha,\, \beta\}\cap\Phi$ is isomorphic to the root system of type $B_2$.
Direct computation yields that $\{\alpha,\, \beta\}=\{e_i,\, -e_i+e_j\}$ for some $i,j$ with $1\leqslant i<j\leqslant n$.
Since $\ell(w| _{\Phi\cap U})=2$, it follows that ${\rm sts}(w_iw_j)=\underline{2}1$ or $2\underline{1}$. So $w$ contains the patterns $\underline{2}1$ or $2\underline{1}$.

For the case (i), there exist $\alpha,\, \beta,\, \gamma\in \Phi^+$ as simple roots (in this order)
such that the subsystem ${\rm span}\{\alpha,\, \beta,\, \gamma\}\cap\Phi$ is isomorphic to the root system of type $A_3$.
So $${\rm span}\{\alpha,\, \beta,\, \gamma\}\cap \Phi^{+}=\{\alpha, \beta,\gamma, \alpha+\beta, \beta+\gamma,\alpha+\beta+\gamma\},$$
from which we see that either
\begin{align}\label{eq:Iphiw11}
I_\Phi(w)\cap V=\{\alpha,\, \alpha+\beta+\gamma,\, \gamma\}\quad \text{and}\quad I_\Phi(w)\cap\{\alpha+\beta,\, \beta,\, \beta+\gamma\}=\emptyset,
\end{align}
or
\begin{align}\label{eq:Iphiw22}
I_\Phi(w)\cap V=\{\alpha+\beta,\, \beta,\, \beta+\gamma\}\quad\, \text{and}\quad I_\Phi(w)\cap\{\alpha,\, \alpha+\beta+\gamma,\, \gamma\}=\emptyset.
\end{align}
{Let
\begin{align*}
\Phi^+_1=\{-e_i+ e_j \mid 1\leq i<j \leq n\}\, , \,\Phi^+_2=\{e_i+ e_j \mid 1\leq i<j \leq n\},~~ \text{and}~~\Phi^+_3=\{e_i\mid i\in[n]\}.
\end{align*}
Then $\Phi^+=\Phi^+_1\uplus\Phi^+_2\uplus\Phi^+_3$.
A simple computation yields that $\{\alpha,\, \beta,\, \gamma\}\cap \Phi^+_3=\emptyset$, and the cardinality of $\{\alpha,\, \beta,\, \gamma\}\cap\Phi^+_1$ is either two or three.

If $\#\{\{\alpha,\, \beta,\, \gamma\}\cap\Phi^+_1\}=2$, then there will be the following three cases.

{\bf Case 1:} $\alpha,\, \beta\in \Phi^+_1$ and $\gamma\in\Phi^+_2$.
Let $\alpha=-e_i+e_j$, $\beta=-e_k+e_l$, $\gamma=e_p+e_q$ for some $i<j$, $k<l$, $p<q$ in $[n]$.
Notice that $\alpha+\beta\in \Phi^{+}$ implies that either $i=l$ or $j=k$, $\beta+\gamma\in \Phi^{+}$ implies that either $k=p$ or $k=q$, and $\alpha+\beta+\gamma\in \Phi^{+}$ implies that $i,k\in\{j,l,p,q\}$.
For the case $j=k=p$, we have $\alpha+\beta+\gamma=-e_i+e_l+e_j+e_q$.
Since $i<j=k<l$ and $\alpha+\beta+\gamma\in \Phi^{+}$, we get $i=q$, contradicting $i<j=p<q$.
Similarly, the case $j=k=q$ is impossible.
For the case $i=l$ and $k=p$ we have
\begin{align*}
\alpha=-e_i+e_j,\ \beta=-e_k+e_i,\ \gamma=e_k+e_q, \ \text{where}\ k<i<j, \ k<q.
\end{align*}
If Eq.\,\eqref{eq:Iphiw11} holds, then
\begin{align*}
w_i>w_j,\ w_j+w_q<0,\ w_k+w_q<0,\ w_k<w_j,\ w_k<w_i,\  \text{and}\ w_i+w_q>0,
\end{align*}
that is, $w_k<w_j<\underline{w_{q}}<w_i$.
Similarly, if  Eq.\,\eqref{eq:Iphiw22} holds, then $w_i<\underline{w_{q}}<w_j<w_k$.
It is routine to verify that $w$ must {contain} $\underline{2}1$ or $2\underline{1}$ in the cases: $k<q<i<j, k<i<q<j$ and $k<i<j<q$.
%, so $w$ contains at least one of the patterns $\underline{2}1$ and $2\underline{1}$.
For the case $i= l$ and $k=q$, a completely analogous argument shows that $w$ also contains at least one of the patterns $\underline{2}1$ and $2\underline{1}$.

{\bf Case 2:} $\beta, \gamma\in \Phi^+_1$ and $\alpha\in\Phi^+_2$. The argument is analogous to that of Case $1$ by interchanging the roles of $\alpha$ and $\gamma$.

{\bf Case 3:} $\alpha, \gamma \in \Phi^+_1$ and $\beta\in \Phi^+_2$.
Let $\alpha=-e_i+e_j$, $\beta=e_k+e_l$, and $\gamma=-e_p+e_q$ where $i<j$, $k<l$, and $p< q$ in $[n]$.
Since $\alpha+\gamma\not\in \Phi^{+}$, we have $i\neq q$ and $j\neq p$, which together with $\alpha+\beta+\gamma\in \Phi^{+}$ implies $e_i+e_p=e_k+e_l$, and hence $\{i,p\}=\{k,l\}$.
Without loss of generality we may assume that $i=k$ and $p=l$. Then
$$
\alpha=-e_i+e_j,\ \beta=e_i+e_p, \ \gamma=-e_p+e_q, \text{where}\ i<j, \ i<p<q,\ j\neq p,q.
$$
Applying Eq.\,\eqref{eq:Iphiw11} we see that  $\underline{w_p}<w_j<\underline{w_{q}}<w_i$, and applying Eq.\,\eqref{eq:Iphiw22} we obtain that $w_i<\underline{w_q}<w_j<\underline{w_{p}}$.
Thus, $w$ contains at least one of the patterns $\underline{2}1$ and $2\underline{1}$.

If $\#\{\{\alpha,\, \beta,\, \gamma\}\cap\Phi^+_1\}=3$, then $\alpha,\, \beta,\, \gamma\in\Phi^+_1$. A simple computation shows that there exist $i,j,k,l$ with $1\leq i< j< k< l\leq n$ such that $\alpha=-e_i+e_j$, $\beta=-e_j+e_k$ and $\gamma=-e_k+e_l$.} Then Eq.\,\eqref{eq:Iphiw11} yields that ${\rm st}(w_iw_jw_kw_l)=3142$, while  Eq.\,\eqref{eq:Iphiw22} yields that ${\rm st}(w_iw_jw_kw_l)=2413$.
From Example \ref{exmp31422413} we see that all of the elements of $\{u\in \mathfrak{B}_4\mid {\rm st}(u)=3142\}\cup\{u\in \mathfrak{B}_4\mid {\rm st}(u)=2413\}$ contain the patterns $2\underline{1}$ or $\underline{2}1$ except for $3142$, $2413$, $\underline{3}\,\underline{1}\,\underline{4}\,\underline{2}$ and $\underline{2}\,\underline{4}\,\underline{1}\,\underline{3}$.
Thus, $w$ contains at least one of the patterns $\underline{2}1$, $2\underline{1}$, $3142$, $2413$, $\underline{3}\,\underline{1}\,\underline{4}\,\underline{2}$, and $\underline{2}\,\underline{4}\,\underline{1}\,\underline{3}$, completing the proof.
\end{proof}

The proof of Theorem \ref{maintheorem} will be reduced to a special class of non-separable elements in $\mathfrak{B}_n$.
An element $w$ of a Weyl group $W$ is called \emph{minimal non-separable} if $w$ is not separable, but $w_J\in W_J$ is separable for all $J\subsetneq\Delta$. By \cite[Proposition 5.2]{GG20aam}, any pattern {which} occurs in a separable element is separable,
so it suffices to require that  $J=\Delta\backslash\{\alpha\}$, where $\alpha$ ranges over all elements of $\Delta$.

\begin{lemma}\label{lem:minonsepBnspec}
If $w=w_1w_2\cdots w_n\in \mathfrak{B}_n$ {is} minimal non-separable, then ${\rm st}(w)$ and ${\rm sts}(w_1w_2 \cdots w_{n-1})$ are separable.	
\end{lemma}
\begin{proof}
Let $J_i=\Delta\backslash\{\alpha_i\}$, where $i\in[0,n-1]$. Then ${\rm st}(w)=w_{J_0}$ and ${\rm sts}(w_1w_2 \cdots w_{n-1})=w_{J_{n-1}}$, so they are separable.
\end{proof}

\begin{lemma}\label{minonsepBn}
	Let $w=w_1w_2\cdots w_n\in \mathfrak{B}_n$. Then $w$ is minimal non-separable if and only if the following conditions  hold:
	\begin{enumerate}
		\item\label{minonsepBn:item1} $w_1w_2 \cdots w_{n-1}$  avoids $\underline{2}1$ and $2\underline{1}$, $w$ avoids  $3142,2413,\underline{3}\,\underline{1}\,\underline{4}\,\underline{2}$, and $\underline{2}\,\underline{4}\,\underline{1}\,\underline{3}$;
		\item\label{minonsepBn:item2} there exists $i\in[n-1]$ such that
		\begin{align*}
			{\rm sts}(w_iw_n)=
			\begin{cases}
				\underline{2}1,&{\rm if}\ w_n>0,\\
				2\underline{1},&{\rm if}\ w_n<0,
			\end{cases}
		\end{align*}
		and for any $j,k,l$ with $1\leq j<k<l<n$, we have
		\begin{align}\label{minonsepBn:item2wijkn}
			{\rm sts}(w_jw_kw_lw_n)\notin
			\begin{cases}
				\{13\underline{4}2,2\underline{3}41,\underline{1}3\underline{4}2,\underline{2}3\underline{4}1\},&{\rm if}\ w_n>0,\\
				\{\underline{1}\,\underline{3}4\underline{2},\underline{2}3\underline{4}\,\underline{1},
				1\underline{3}4\underline{2},2\underline{3}4\underline{1}\},&{\rm if}\ w_n<0.
			\end{cases}
		\end{align}
	\end{enumerate}
\end{lemma}
\begin{proof}
Assume that $w$ is a minimal non-separable element of $\mathfrak{B}_n$.
Then, by Lemma \ref{lem:minonsepBnspec}, ${\rm st}(w)$ and
${\rm sts}(w_1w_2\cdots w_{n-1})$ are separable. Hence $w_1w_2\cdots w_{n-1}$ avoids the patterns $\underline{2}1$ and $2\underline{1}$ by
Lemma \ref{sepBnavoidpattern}, and
${\rm st}(w)$ avoids the patterns $3142$ and $2413$.
This implies that $w$ avoids all patterns from the sets
$\{u\in \mathfrak{B}_4\mid {\rm st}(u)=3142\}$ and $\{u\in \mathfrak{B}_4\mid {\rm st}(u)=2413\}$.
In particular, it follows from Example \ref{exmp31422413} that $w$ avoids the patterns $3142$, $2413$, $\underline{3}\,\underline{1}\,\underline{4}\,\underline{2}$ and $\underline{2}\,\underline{4}\,\underline{1}\,\underline{3}$, as well as the eight patterns on the right-hand side of Eq.\,\eqref{minonsepBn:item2wijkn}.
Thus, both part \eqref{minonsepBn:item1} and Eq.\,\eqref{minonsepBn:item2wijkn} are true.

Note that $w$ is non-separable. It follows from Lemma \ref{sepBnavoidpattern} that $w$ contains at least one of the patterns $\underline{2}1$, $2\underline{1}$, $3142$, $2413$, $\underline{3}\,\underline{1}\,\underline{4}\,\underline{2}$, and $\underline{2}\,\underline{4}\,\underline{1}\,\underline{3}$.
Combining this with  part \eqref{minonsepBn:item1} implies that
$w$ must contain the patterns $\underline{2}1$ or $2\underline{1}$, that is, there exist $i$ and $j$ with $1\leq i<j\leq n$ such that ${\rm sts}(w_iw_j)=\underline{2}1$ or $2\underline{1}$.
But $w_1w_2\cdots w_{n-1}$ avoids the patterns $\underline{2}1$ and $2\underline{1}$, so $j=n$, and hence ${\rm sts}(w_iw_n)=\underline{2}1$ if $w_n>0$,
and ${\rm sts}(w_iw_n)=2\underline{1}$ otherwise. Thus, the proof of part \eqref{minonsepBn:item2} follows.

Conversely, according to part \eqref{minonsepBn:item2}, $w$ contains the patterns $\underline{2}1$ or $2\underline{1}$ so that $w$ is a non-separable element.
Let $J_i=\Delta\backslash\{\alpha_i\}$, where $i\in[0,n-1]$.
It remains to show that $w_{J_i}$ is separable in $\mathfrak{B}_{J_i}$ for all $i\in [0,n-1]$.
For each $i\in [0,n-1]$, we have
$$w_{J_i}={\rm sts}(w_1\cdots w_i)\times {\rm st}(w_{i+1}\cdots w_n)\in \mathfrak{B}_{i}\times\mathfrak{S}_{n-i},$$
where $\mathfrak{B}_{0}=\mathfrak{S}_{0}=\{e\}$ is the trivial group.
From part \eqref{minonsepBn:item1} we see that ${\rm sts}(w_1\cdots w_i)$ avoids the patterns $\underline{2}1$, $2\underline{1}$, $3142$, $2413$, $\underline{3}\,\underline{1}\,\underline{4}\,\underline{2}$, and $\underline{2}\,\underline{4}\,\underline{1}\,\underline{3}$, and hence it is separable according to
Lemma \ref{sepBnavoidpattern}. To prove $w_{J_i}$ is separable, it is enough to show that ${\rm st}(w_{i+1}\cdots w_n)$ is separable by Definition \ref{separabledef}.
This is equivalent to the statement that $w_{i+1}\cdots w_n$ avoids all patterns $u\in \mathfrak{B}_4$ that satisfy ${\rm st}(u)=3142$ or $2413$.
It follows from part \eqref{minonsepBn:item1} that $w_{i+1}\cdots w_{n-1}$  avoids the patterns $\underline{2}1$ and $2\underline{1}$, $w_{i+1}\cdots w_{n}$  avoids the patterns $3142,2413,\underline{3}\,\underline{1}\,\underline{4}\,\underline{2}$, and $\underline{2}\,\underline{4}\,\underline{1}\,\underline{3}$.
Compare these with Example \ref{exmp31422413},
it will be clear that we only need to show $w_{i+1}\cdots w_{n}$ does not contain the eight patterns on the right-hand side of Eq.\,\eqref{minonsepBn:item2wijkn}.
Suppose to the contrary that $w_{i+1}\cdots w_{n}$ contains one of the eight patterns, say $v=v_1v_2v_3v_4$.
Then $v$ contains at least one of the patterns $\underline{2}1$ and $2\underline{1}$. Thus, part \eqref{minonsepBn:item1} guarantees that any occurrence of the pattern $v$ in $w_{i+1}\cdots w_{n}$ must use $w_n$.
But this contradicts Eq.\,\eqref{minonsepBn:item2wijkn}. So ${\rm st}(w_{i+1}\cdots w_n)$
is separable, completing the proof.
\end{proof}

According \cite[Corollary 2]{GG20am}, if an element $w$ of a Weyl group $W$ is separable, then $w_0w$, $ww_0$, and $w^{-1}$ are also separable, where $w_0$ is the longest element of $W$.
Since the longest element $w_0$ of $\mathfrak{B}_n$ is $\underline{1}\,\underline{2}\cdots \underline{n}$, it follows from Lemma \ref{minonsepBn} that
$w\in \mathfrak{B}_n$ is minimal non-separable implies $w_0w$ and $ww_0$ are also minimal non-separable.
So we get the following result.

\begin{lemma}\label{ww0w-1separableprop}
Let $w\in \mathfrak{B}_n$ be non-separable (respectively, minimal non-separable). Then $w_0w$ and $ww_0$ are also non-separable (respectively, minimal non-separable).
\end{lemma}

For a minimal non-separable element $w$ of $\mathfrak{B}_n$,  according to Lemma \ref{sepBnavoidpattern}, its inverse $w^{-1}$ is also non-separable.
However, $w^{-1}$ may be not minimal in general.
For example, by Lemma \ref{minonsepBn}, the element $w=\underline{2}3451$ is minimal non-separable. But $w^{-1}=5\underline{1}234$ is not minimal non-separable,
since ${\rm sts}(w^{-1}_1w^{-1}_2)=2\underline{1}$.
The following lemma describes the structure of the minimal non-separable signed permutations whose inverse elements are also minimal non-separable.

\begin{lemma}\label{w-1}
Let $w=w_1w_2\cdots w_n\in \mathfrak{B}_n$ be minimal non-separable. Then $w^{-1}$ is minimal non-separable if and only if the following conditions hold:
\begin{enumerate}
  \item\label{w-1:item1} there exists $i\in[n-1]$ such that $|w_i|=n$ and ${\rm sts}(w_1\cdots \hat{w_i}\cdots w_{n})\in \mathfrak{B}_{n-1}$ is separable, where the symbol $\hat{w_i}$ indicates that $w_i$ has been omitted;
  \item\label{w-1:item2} for any $j,k,l$ with $1\leq j<k<l<n$, we have ${\rm sts}(w_jw_kw_lw_n)\notin \{\underline{1}\,\underline{4}\,\underline{2}3,1\underline{4}\,\underline{2}3\}$ if $w_n>0$,
     and ${\rm sts}(w_jw_kw_lw_n)\notin \{142\underline{3},\underline{1}42\underline{3}\}$ if $w_n<0$.
\end{enumerate}
\end{lemma}
\begin{proof}
Suppose $w^{-1}$ is minimal non-separable. It follows from  Lemma \ref{lem:minonsepBnspec} that ${\rm st}(w^{-1})$ and
${\rm sts}(w^{-1}_1w^{-1}_2\cdots w^{-1}_{n-1})$ are separable. But ${\rm sts}(w_1\cdots \hat{w_i}\cdots w_n)=({\rm sts}(w^{-1}_1w^{-1}_2\cdots w^{-1}_{n-1}))^{-1}$, so it is separable by Lemma \ref{sepBnavoidpattern}, where $|w_i|=n$.
Since $w$ is minimal non-separable, it follows from Lemma \ref{minonsepBn}\eqref{minonsepBn:item2} that $|w_n|\neq n$, so that $i\in[n-1]$.
Thus part \eqref{w-1:item1} is true.
Assume that there exist $j,\,k,\,l$ with $1\leq j<k<l<n$  such that
${\rm sts}(w_jw_kw_lw_n)\in \{\underline{1}\,\underline{4}\,\underline{2}3,1\underline{4}\,\underline{2}3,142\underline{3},\underline{1}42\underline{3}\}$.
Let  $p=|w_j|$, $q=|w_l|$, $r=|w_n|$, and $s=|w_k|$. Then $ p<q<r<s$, and
$${\rm sts}(w^{-1}_pw^{-1}_qw^{-1}_rw^{-1}_s)\in \{(\underline{1}\,\underline{4}\,\underline{2}3)^{-1},(1\underline{4}\,\underline{2}3)^{-1},(142\underline{3})^{-1},
(\underline{1}42\underline{3})^{-1}\}=\{\underline{1}\,\underline{3}4\underline{2},1\underline{3}4\underline{2},13\underline{4}2,\underline{1}3\underline{4}2\},$$
so that ${\rm st}(w^{-1}_pw^{-1}_qw^{-1}_rw^{-1}_s)$ is either
$3142$ or $2413$, contradicting the fact that ${\rm st}(w^{-1})$ is separable. Hence part \eqref{w-1:item2} holds.

Conversely, suppose parts \eqref{w-1:item1} and \eqref{w-1:item2} hold. We assume that $w_n>0$; the case $w_n<0$ is similar.
By Lemma \ref{sepBnavoidpattern} and part \eqref{w-1:item1},
${\rm sts}(w^{-1}_1w^{-1}_2\cdots w^{-1}_{n-1})=({\rm sts}(w_1\cdots \hat{w_i}\cdots w_n))^{-1}$
is separable. Therefore, $w^{-1}_1w^{-1}_2\cdots w^{-1}_{n-1}$ avoids the patterns $\underline{2}1$ and $2\underline{1}$ by Lemma \ref{sepBnavoidpattern}.
Since $w$ is minimal non-separable, according to Lemma \ref{minonsepBn}\eqref{minonsepBn:item1}, $w$ avoids the patterns $3142$, $2413$, $\underline{3}\,\underline{1}\,\underline{4}\,\underline{2}$, and $\underline{2}\,\underline{4}\,\underline{1}\,\underline{3}$, and hence so does $w^{-1}$.

Since we have already supposed $w_n>0$, it follows from Lemma \ref{minonsepBn}\eqref{minonsepBn:item2} that
${\rm sts}(w_jw_n)=\underline{2}1$ for some
$j\in[n-1]$.
Then ${\rm sts}(w^{-1}_pw^{-1}_q)=2\underline{1}$, where $p=w_n$ and $q=|w_j|$.
But $w^{-1}_1w^{-1}_2\cdots w^{-1}_{n-1}$ avoids the pattern $2\underline{1}$, so $q=n$, and hence ${\rm sts}(w^{-1}_pw^{-1}_n)=2\underline{1}$. In particular, we have $w_n^{-1}<0$.
Assume that there exist $j,k,l$ with $1\leq j<k<l<n$ such that
${\rm sts}(w^{-1}_jw^{-1}_kw^{-1}_lw^{-1}_n)\in\{\underline{1}\,\underline{3}4\underline{2},\underline{2}3\underline{4}\,\underline{1},
				1\underline{3}4\underline{2},2\underline{3}4\underline{1}\}$.
If ${\rm sts}(w^{-1}_jw^{-1}_kw^{-1}_lw^{-1}_n)\in\{\underline{1}\,\underline{3}4\underline{2},1\underline{3}4\underline{2}\}$, then
${\rm sts}(w_{r_1}w_{r_2}w_{r_3}w_{r_4})\in\{\underline{1}\,\underline{4}\,\underline{2}3,1\underline{4}\,\underline{2}3\}$,
where $r_1=|w^{-1}_j|<r_2=|w^{-1}_n|<r_3=|w^{-1}_k|<r_4=|w^{-1}_l|$.
It follows from  part \eqref{w-1:item2} that $r_4<n$.
Then $w_1w_2\cdots w_{n-1}$ contains the pattern ${\rm sts}(w_{r_2}w_{r_4})=\underline{2}1$.
If ${\rm sts}(w^{-1}_jw^{-1}_kw^{-1}_lw^{-1}_n)\in\{\underline{2}3\underline{4}\,\underline{1},2\underline{3}4\underline{1}\}$, then
${\rm sts}(w_{r_1}w_{r_2}w_{r_3}w_{r_4})\in\{\underline{4}\,\underline{1}2\,\underline{3},\underline{4}1\underline{2}3\}$,
where $r_1=|w^{-1}_n|<r_2=|w^{-1}_j|<r_3=|w^{-1}_k|<r_4=|w^{-1}_l|$. So $w_{r_1}w_{r_2}w_{r_3}$
contains the pattern $\underline{2}1$, and hence $w_1w_2\cdots w_{n-1}$ contains $\underline{2}1$.
Thus in either case we have $w_1w_2\cdots w_{n-1}$ contains $\underline{2}1$, contradicting the fact that $w$ is
minimal non-separable by Lemma \ref{minonsepBn}\eqref{minonsepBn:item1}. This proves that ${\rm sts}(w^{-1}_jw^{-1}_kw^{-1}_lw^{-1}_n)\notin \{\underline{1}\,\underline{3}4\underline{2},\underline{2}3\underline{4}\,\underline{1}, 1\underline{3}4\underline{2},2\underline{3}4\underline{1}\}$ for any $j,k,l$ with $1\leq j<k<l<n$. Therefore, $w^{-1}$ is minimal non-separable by Lemma \ref{minonsepBn}, and we are done.
\end{proof}

\section{Splittings in Weyl groups of type $B$}\label{sec:Proofofthemainythm}

This section is devoted to proving Theorem \ref{maintheorem}, by means of Lemmas \ref{notranksymmetric-(n-1)} and \ref{notranksymmetric(n-1)} below.
One of the crucial steps of our proof is to show that $\Lambda^{L}_{w}$  is not rank-symmetric when $w$ and $w^{-1}$  are both minimal non-separable with some additional conditions.
For this we need to analyze the detailed structure of  minimal non-separable elements of $\mathfrak{B}_n$.

\begin{lemma}\label{nwn-1nn-1n}
Let $w=w_1w_2\cdots w_n\in \mathfrak{B}_n$. If $w$ and $w^{-1}$ are both minimal non-separable, then either $|w_n|=n-1$ or $|w_{n-1}|=n$.
\end{lemma}
\begin{proof}
Since $w$ is minimal non-separable, from Lemma \ref{lem:minonsepBnspec} we see that ${\rm sts}(w_1w_2\cdots w_{n-1})$ is separable, and by Lemma \ref{minonsepBn}\eqref{minonsepBn:item2}
there exists $i'\in[n-1]$ such that ${\rm sts}(w_{i'}w_n)\in\{\underline{2}1,2\underline{1}\}$.
Since $w^{-1}$ is also minimal non-separable, we deduce from Lemma \ref{w-1} that there exists $i\in[n-1]$ such that
$|w_i|=n$ and  ${\rm sts}(w_1\cdots \hat{w_i}\cdots w_n)$ is separable.
According to Lemma \ref{sepBnavoidpattern}, ${\rm sts}(w_1\cdots \hat{w_i}\cdots w_n)$ avoids the patterns $\underline{2}1$ and $2\underline{1}$, so
$i$ is the unique index such that ${\rm sts}(w_iw_n)\in\{\underline{2}1,2\underline{1}\}$.

Assume that $|w_n|\neq n-1$. We need to show that $|w_{n-1}|=n$, that is, $i=n-1$. Suppose to the contrary that $i\neq n-1$.
Note that ${\rm sts}(w_iw_n)\in\{\underline{2}1,2\underline{1}\}$ guarantees that $|w_{n}|\neq n$, so $i<n-1$.
Let $j\in[n-1]\backslash\{i\}$ be such that $|w_j|=n-1$.  By the uniqueness of $i$, we see that $w_{j}$ and $w_n$ have the same signs, so that $w_{j}$ and $w_i$ have the opposite signs.
If $i<j$, then
${\rm sts}(w_iw_j)\in\{\underline{2}1,2\underline{1}\}$,
if $j<i$, then either ${\rm sts}(w_iw_{n-1})\in\{\underline{2}1,2\underline{1}\}$ or ${\rm sts}(w_jw_{n-1})\in\{\underline{2}1,2\underline{1}\}$.
Thus, $w_1w_2 \cdots w_{n-1}$ must contain the patterns $\underline{2}1$ or $2\underline{1}$, contradicting
Lemma \ref{minonsepBn}\eqref{minonsepBn:item1}, proving that $i=n-1$.
\end{proof}

\begin{lemma}\label{lem:wjk>0<0uponwn+-}
Let $w=w_1w_2\cdots w_n\in \mathfrak{B}_n$ with $|w_n|=n-1$ and $|w_i|=n$ for some $i\in[n-2]$, and suppose that $w$ and $w^{-1}$ are both minimal non-separable. Then for any $j$ and $k$ with $1\leq j<i<k<n$, we have
$w_j>w_k>0$ if $w_n=\underline{n-1}$, and $w_j<w_k<0$ if $w_n=n-1$.
\end{lemma}
\begin{proof}
Let $j$ and $k$ be such that $1\leq j<i<k<n$.
If $w_n=\underline{n-1}$, then, by Lemma \ref{minonsepBn}\eqref{minonsepBn:item2}, $w_{i}=n$.
But Lemma \ref{minonsepBn}\eqref{minonsepBn:item1} says that ${\rm sts}(w_iw_k)\neq 2\underline{1}$, so $w_{k}>0$.
Because $w^{-1}$ is also minimal non-separable,  it follows from Lemma \ref{w-1}\eqref{w-1:item2} that ${\rm sts}(w_jw_iw_kw_n)\notin \{142\underline{3},\underline{1}42\underline{3}\}$,
so that $|w_j|>w_k>0$. By Lemma \ref{minonsepBn}\eqref{minonsepBn:item1}, ${\rm sts}(w_jw_k)\neq \underline{2}1$, and hence $w_j>w_k>0$.
Substitute $w_0w$ for $w$, a completely analogous argument shows that $w_j<w_k<0$ if $w_n=n-1$, as claimed.
\end{proof}

We now show that $\Lambda^{L}_{w}$  is not rank-symmetric if $w$ and $w^{-1}$ are both minimal non-separable with $w_n=\underline{n-1}$ and $w_i=n$ for some $i\in [n-2]$.

\begin{lemma}\label{notranksymmetric-(n-1)}
Let $w=w_1w_2\cdots w_n\in \mathfrak{B}_n$ with $w_n=\underline{n-1}$ and $w_i=n$ for some $i\in [n-2]$, and let $f(q)=\Lambda_w^L(q)$. If $w$ and $w^{-1}$ are both minimal non-separable, then
\begin{align*}
[q^d]f=\begin{cases}
[q^{\ell(w)-d}]f, &{\rm if}\  0\leq d\leq i-1,\\
[q^{\ell(w)-d}]f+1, &{\rm if}\ d=i.
\end{cases}
\end{align*}
In particular, $\Lambda^{L}_{w}$  is not rank-symmetric.
\end{lemma}
\begin{proof}
Since $w_n=\underline{n-1}$ and $w_i=n$, we have
$w=w_1\cdots w_{i-1}nw_{i+1}\cdots w_{n-1}\underline{(n-1)}$.
It follows from Lemma \ref{lem:wjk>0<0uponwn+-} that
\begin{align*}
\{w_1, \ldots, w_{i-1}\}=\{n-i, \ldots, n-2\}~\text{and}~ \{w_{i+1}, \ldots, w_{n-1}\}=\{1, \ldots, n-i-1\},
\end{align*}
so that ${\rm Neg}(w)=\{n\}$, ${\rm Nsp}(w)=\{(j,n)\mid j\in[n-1]\backslash\{i\}\}$, and
\begin{align*}%\label{eq:invw}
{\rm Inv}(w)={\rm Inv}(w_1\cdots w_{i-1})&\uplus {\rm Inv}(w_{i+1}\cdots w_{n-1})
\uplus \{(j,k)\mid 1\leq j\leq i<k<n\}\uplus\{(l,n)\mid l\in[n-1]\}.
\end{align*}
Let $w'={\rm st}(w)$. Then
$w'=w'_1\cdots w'_{i-1}nw'_{i+1}\cdots w'_{n-1}1$,
where $w'_j=w_j+1$ for $j\in[n-1]\backslash\{i\}$. In particular, we have
\begin{align}\label{eq:w'1inw'shap}
\{w'_1, \ldots, w'_{i-1}\}=\{n-i+1, \ldots, n-1\}~\text{and}~ \{w'_{i+1}, \ldots, w'_{n-1}\}=\{2, \ldots, n-i\}.
\end{align}
Thus, ${\rm Inv}(w')={\rm Inv}(w)$, and hence $\ell(w)=\ell(w')+n-1$ by Eq.\,\eqref{eq:eqell(w)}.

Let $g(q)=[e,w']_L(q)$. By Lemma \ref{lem:minonsepBnspec}, $w'$ is separable, which together with Lemma \ref{rank-symmetric} yields that $g(q)$ is symmetric and unimodal.
For any $u\in[e,w]_L$ with $\ell(u)\leq i$, the set ${\rm Neg}(u)$ must be empty. Otherwise, ${\rm Neg}(u)=\{n\}$ so that
$$\ell(u)\geq \#{\rm Inv}(u)+\#{\rm Neg}(u)>\#{\rm Inv}(u)\geq\#\{(j,n)|j\in[n-1]\}=n-1>i,$$
a contradiction, proving ${\rm Neg}(u)=\emptyset$. Thus, ${\rm Nsp}(u)=\emptyset$ and
${\rm Inv}(u)\subseteq {\rm Inv}(w)={\rm Inv}(w')$, which together with Lemma \ref{leftweakorder} yields that $u\leq_Lw'$.
Since $w'\leq_L w$, we see that for any $u\in \mathfrak{B}_n$ with $\ell(u)\leq i$, $u\leq_Lw$ is equivalent to $u\leq_L w'$, and hence for any $d\leq i$, we have
\begin{align*}
[q^d]f=[q^d]g=[q^{\ell(w')-d}]g.
\end{align*}
Thus it suffices to show that
\begin{equation*}%\label{eq:qw'-dgf}
%\begin{split}
[q^{\ell(w')-d}]g=\begin{cases}
[q^{\ell(w)-d}]f, &{\rm if}\  0\leq d\leq i-1,\\
[q^{\ell(w)-d}]f+1, &{\rm if}\ d=i.
\end{cases}
%\end{split}
\end{equation*}

Let $v= s_{n-1}s_{n-2}\cdots s_{n-i}w'$, and let
\begin{align*}
A=\{ x\in[e,w']_L\mid \ell(x)\geq\ell(w')-i\}\quad \text{and} \quad
B=\{x\in[e,w]_L\mid \ell(x)\geq\ell(w)-i\}.
\end{align*}
It follows from Eq.\,\eqref{eq:w'1inw'shap} that
$${\rm Inv}(w')={\rm Inv}(v)\uplus \{(k_1,k_2)\mid k_1\in[i], w_{k_2}'=n-i\},$$
and hence $$\ell(w')-\ell(v)=i=\ell(s_{n-i}\cdots s_{n-2}s_{n-1})=\ell(w'v^{-1}).$$
So $v\in A$. Define a map
\begin{align*}
	\psi: A\backslash\{v\} \rightarrow B,\quad x\mapsto x',
\end{align*}
where $x'\in \mathfrak{B}_n$ such that $x'_n=\underline{n-x_n}$ and ${\rm sts}(x'_1x'_2\cdots x'_{n-1})={\rm sts}(x_1x_2\cdots x_{n-1})$.
See Figure $1$ for an example in $\mathfrak{B}_5$, where $w=2351\underline{4}$.
It thus suffices to show that $\psi$ is bijective and that $\ell(x')-\ell(x)=\ell(w)-\ell(w')=n-1$ for any $x\in A\backslash\{v\}$.

\begin{figure}[t]
\begin{center}
\begin{tikzpicture}\label{Figure 2.}
\node (2351underline{4}) at (0,6.6) {{\tiny$2\thinspace3\thinspace5\thinspace1\thinspace \underline{4}$}};

\node (1352ubderline{4}) at (-1,5.4) {{\tiny$1\thinspace3\thinspace5\thinspace2\thinspace \underline{4}$}};
\node (2451underline{3}) at (1,5.4) {{\tiny$2\thinspace4\thinspace5\thinspace1\thinspace \underline{3}$}};

\node (1253underline{4}) at (-2,4.2) {{\tiny$1\thinspace2\thinspace5\thinspace3\thinspace \underline{4}$}};
\node (1452underline{3}) at (0,4.2) {{\tiny$1\thinspace4\thinspace5\thinspace2\thinspace \underline{3}$}};
\node (3451underline{2}) at (2,4.2) {{\tiny$3\thinspace4\thinspace5\thinspace1\thinspace \underline{2}$}};

\node (1254underline{3}) at (-2,3) {{\tiny$1\thinspace2\thinspace5\thinspace4\thinspace \underline{3}$}};
\node (1453underline{2}) at (0,3) {{\tiny$1\thinspace4\thinspace5\thinspace3\thinspace \underline{2}$}};
\node (3452underline{1}) at (2,3) {{\tiny$3\thinspace4\thinspace5\thinspace2\thinspace \underline{1}$}};

\node (1245underline{3}) at (-3,1.8) {{\tiny$1\thinspace2\thinspace4\thinspace5\thinspace \underline{3}$}};
\node (1354underline{2}) at (-1,1.8) {{\tiny$1\thinspace3\thinspace5\thinspace4\thinspace \underline{2}$}};
\node (2453underline{1}) at (1,1.8) {{\tiny$2\thinspace4\thinspace5\thinspace3\thinspace \underline{1}$}};
\node (34521) at (3,1.8) {{\tiny$3\thinspace4\thinspace5\thinspace2\thinspace 1$}};

\node (1345underline{2}) at (-3,0.6) {{\tiny$1\thinspace3\thinspace4\thinspace5\thinspace \underline{2}$}};
\node (2354underline{1}) at (-1,0.6) {{\tiny$2\thinspace3\thinspace5\thinspace4\thinspace \underline{1}$}};
\node (24531) at (1,0.6) {{\tiny$2\thinspace4\thinspace5\thinspace3\thinspace 1$}};
\node (34512) at (3,0.6) {{\tiny$3\thinspace4\thinspace5\thinspace1\thinspace 2$}};

\node (2345underline{1}) at (-3,-0.6) {{\tiny$2\thinspace3\thinspace4\thinspace5\thinspace \underline{1}$}};
\node (23541) at (-1,-0.6) {{\tiny$2\thinspace3\thinspace5\thinspace4\thinspace 1$}};
\node (14532) at (1,-0.6) {{\tiny$1\thinspace4\thinspace5\thinspace3\thinspace 2$}};
\node (24513) at (3,-0.6) {{\tiny$2\thinspace4\thinspace5\thinspace1\thinspace 3$}};

\node (23451) at (-3,-1.8) {{\tiny$2\thinspace3\thinspace4\thinspace5\thinspace 1$}};
\node (13542) at (-1,-1.8) {{\tiny$1\thinspace3\thinspace5\thinspace4\thinspace 2$}};
\node (14523) at (1,-1.8) {{\tiny$1\thinspace4\thinspace5\thinspace2\thinspace 3$}};
\node (23514) at (3,-1.8) {{\tiny$2\thinspace3\thinspace5\thinspace1\thinspace 4$}};

\node (13452) at (-3,-3) {{\tiny$1\thinspace3\thinspace4\thinspace5\thinspace 2$}};
\node (12543) at (-1,-3) {{\tiny$1\thinspace2\thinspace5\thinspace4\thinspace 3$}};
\node (13524) at (1,-3) {{\tiny$1\thinspace3\thinspace5\thinspace2\thinspace 4$}};
\node (23415) at (3,-3) {{\tiny$2\thinspace3\thinspace4\thinspace1\thinspace 5$}};

\node (12453) at (-2,-4.2) {{\tiny$1\thinspace2\thinspace4\thinspace5\thinspace 3$}};
\node (12534) at (0,-4.2) {{\tiny$1\thinspace2\thinspace5\thinspace3\thinspace 4$}};
\node (13425) at (2,-4.2) {{\tiny$1\thinspace3\thinspace4\thinspace2\thinspace 5$}};

\node (12354) at (-1,-5.4) {{\tiny$1\thinspace2\thinspace3\thinspace5\thinspace 4$}};
\node (12435) at (1,-5.4) {{\tiny$1\thinspace2\thinspace4\thinspace3\thinspace 5$}};

\node (12345) at (0,-6.6) {{\tiny$1\thinspace2\thinspace3\thinspace4\thinspace 5$}};
\draw[dashed,line width=0.5pt] (-2.8,7.1)--(2.8,7.1)--(2.8,2.55)--(-2.8,2.55)--(-2.8,7.1);
\node (B) at (3.8,4.5) {{\tiny$B$}};
\draw[solid,line width=0.8pt] (2.8,4.5)--(B);

\draw[dashed,line width=0.5pt] (3.8,2.3)--(3.8,-2.3)--(-1.8,-2.3)--(-1.8,-0.5)--(2.7,2.3)--(3.8,2.3);
\node (A) at (5.1,0) {{\tiny$A\backslash\{v\}$}};
\draw[solid,line width=0.8pt] (3.8,0)--(A);

\draw[dashed,line width=0.5pt] (-3.8,-2.3)--(-3.8,-1.3)--(-2.2,-1.3)--(-2.2,-2.3)--(-3.8,-2.3);
\node (v) at (-4.8,-1.8) {{\tiny$v$}};
\draw[solid,line width=0.8pt] (-3.8,-1.8)--(v);

\draw[solid,line width=0.8pt] (2351underline{4})--(1352ubderline{4})--(1253underline{4})- -(1254underline{3})--(1245underline{3})--(1345underline{2})--(2345underline{1})--(23451)--(13452)--(12453)--(12354)--(12345)--(12435)--(13425)--(23415)--(23514)--
(24513)--(34512)--(34521)--(3452underline{1})--(3451underline{2})--(2451underline{3})--(2351underline{4});
\draw[solid,line width=0.8pt] (1352ubderline{4})--(1452underline{3})--(1453underline{2})--(1354underline{2})--(1345underline{2});
\draw[solid,line width=0.8pt] (1254underline{3})--(1354underline{2});
\draw[solid,line width=0.8pt] (2451underline{3})--(1452underline{3});
\draw[solid,line width=0.8pt] (1354underline{2})--(2354underline{1})--(23541)--(13542)--(12543)--(12534)--(13524)--(14523)--(14532)--(24531)--(2453underline{1})--(1453underline{2});
\draw[solid,line width=0.8pt] (2345underline{1})--(2354underline{1})--(2453underline{1})--(3452underline{1});
\draw[solid,line width=0.8pt] (23451)--(23541)--(24531)--(34521);
\draw[solid,line width=0.8pt] (13452)--(13542)--(14532);
\draw[solid,line width=0.8pt] (12534)--(13524)--(23514);
\draw[solid,line width=0.8pt] (14523)--(24513);
\draw[solid,line width=0.8pt] (12453)--(12543);
\draw[solid,line width=0.8pt] (12534)--(12435);
\draw[solid,line width=0.8pt] (13524)--(13425);
\node(Figure 1.)   at (0,-7.6) {Figure $1.$ The left weak order of $[e,2351\underline{4}]_L$.};
\end{tikzpicture}
\end{center}
\end{figure}

Let $x\in A\backslash\{v\}$, and let $x_i=m$. Since ${\rm Inv}(x)\subseteq {\rm Inv}(w')$ and $\{(j,i)\mid j\in[i-1]\}\cap {\rm Inv}(w')=\emptyset$, we have
$[m+1,n]\subseteq\{x_{i+1},\ldots,x_n\}$, so that
$\{(j,k)\mid j\in[i], x_k\in[m+1,n]\}\subseteq {\rm Inv}(w')\setminus{\rm Inv}(x)$. Thus,
$$
i\geq \ell(w')-\ell(x)\geq \#\{(j,k)\mid j\in[i], \, x_k\in[m+1,n]\}=i(n-m),
$$
from which we see that $m\in\{n-1,n\}$. If $m=n-1$, then $\ell(w')-\ell(x)=i$.
Let $k=x^{-1}_n$. Then $x_k=n$ and $k>i$.
Note that $\{(j,k)\mid j\in[i]\}$ is a subset of ${\rm Inv}(w')\setminus{\rm Inv}(x)$, and it contains exactly $i$ elements.
So ${\rm Inv}(w')={\rm Inv}(x)\uplus \{(j,k)\mid j\in[i]\}$, from which we see that
$w_{k}'=n-i$ and ${\rm Inv}(x_1\cdots \hat{x}_k\cdots x_n)={\rm Inv}(w'_1\cdots \hat{w'}_k\cdots w'_n)$.
So ${\rm Inv}(x)={\rm Inv}(v)$ and hence $x=v$, a contradiction. Thus, $m=n$, that is, $x_i=n$, so that $x_i'=n>n-x_n$.
Hence $x_n'=\underline{n-x_n}$ is well defined.
It follows from ${\rm Neg}(w')=\emptyset$ that ${\rm Neg}(x)=\emptyset$, and thus
\begin{align}\label{eq:negx'nspx'}
{\rm Neg}(x')=\{n\}={\rm Neg}(w)\quad {\rm and}\quad {\rm Nsp}(x')=\{(j,n)\mid j\in[n-1]\backslash\{i\}, x_j'<n-x_n\}\subseteq {\rm Nsp}(w).
\end{align}
Since $x\leq_Lw'\leq_Lw$, we see from the definition of $x'$ that
\begin{align}\label{eq:invx'invw}
{\rm Inv}(x')&={\rm Inv}(x_1x_2\cdots x_{n-1})\uplus\{(j,n) \mid j\in[n-1]\}\\
&\subseteq {\rm Inv}(w_1w_2\cdots w_{n-1})\uplus\{(j,n) \mid j\in[n-1]\}={\rm Inv}(w).\notag
\end{align}
This shows $x'\leq _Lw$.
Note that ${\rm Neg}(x)$ and ${\rm Nsp}(x)$ are both empty sets, and that
\begin{align}\label{eq:invxxj>xn}
{\rm Inv}(x)={\rm Inv}(x_1x_2\cdots x_{n-1})\uplus \{(j,n)\mid j\in[n-1], x_j>x_n\}.
\end{align}
Comparing Eqs.\,\eqref{eq:negx'nspx'}, \eqref{eq:invx'invw} and \eqref{eq:invxxj>xn} yields
\begin{align*}
\ell(x')-\ell(x)=1+(n-x_n-1)+(x_n-1)=n-1,
\end{align*}
so that
\begin{align*}
\ell(x')=\ell(x)+n-1\geq \ell(w')-i+n-1=\ell(w)-i.
\end{align*}
Thus, $x'\in B$ and hence $\psi$  is well-defined.

For any distinct elements $x,y$ in $A\backslash\{v\}$, if $x_n\neq y_n$, then $x_n'=\underline{n-x_n}\neq \underline{n-y_n}= y_n'$ so that $x'\neq y'$; if $x_n=y_n$, then
$${\rm Inv}(x_1'x_2'\cdots x_{n-1}')={\rm Inv}(x_1x_2\cdots x_{n-1})\neq {\rm Inv}(y_1y_2\cdots y_{n-1})={\rm Inv}(y_1'y_2'\cdots y_{n-1}')$$ so that ${\rm Inv}(x')\neq {\rm Inv}(y')$. Thus, we always have $x'\neq y'$, that is, $\psi$ is injective.

Let $z'\in B$. If $z'_n>0$, then ${\rm Neg}(z')$ and ${\rm Nsp}(z')$ are empty, so that
\begin{align*}
\ell(w)-\ell(z')\geq \#{\rm Neg}(w)+\#{\rm Nsp}(w)=n-1>i,
\end{align*}
a contradiction, so $z'_n<0$. But $(i,n)\not\in {\rm Nsp}(w)$ implies $(i,n)\not\in {\rm Nsp}(z')$, and hence $\underline{n}<z'_n<0$.
Define $z\in \mathfrak{B}_n$ such that $z_n=n+z_n'$ and ${\rm sts}(z_1z_2\cdots z_{n-1})={\rm sts}(z_1'z_2'\cdots z_{n-1}')$.
Then ${\rm Neg}(z)$ and ${\rm Nsp}(z)$ are empty, and
$$
{\rm Inv}(z)={\rm Inv}(z_1'z_2'\cdots z_{n-1}')\uplus \{(j,n)\mid j\in[n-1], z_j>z_n\}\subseteq{\rm Inv}(w').
$$
Hence $z\leq_Lw'$. It follows from Eq.\,\eqref{eq:eqell(w)} that
$\ell(z)=\#{\rm Inv}(z)=\#{\rm Inv}(z_1'z_2'\cdots z_{n-1}')+n-z_n$
and
$$\ell(z')={\rm Inv}(z_1'z_2'\cdots z_{n-1}')+\# \{(j,n)\mid j\in[n-1]\}+n-z_n={\rm Inv}(z_1'z_2'\cdots z_{n-1}')+2n-z_n-1.$$
Thus, $\ell(z)=\ell(z')-n+1\geq \ell(w)-i-n+1=\ell(w')-i$,
proving that $z\in A$.
We claim that $z_i=n$, from which we get $z\neq v$ since $v_i=n-1$.
Let $z_{k}'=n$ for some $k\in[n]$.
Since $\{(j,i)\mid 1\leq j<i\}\cap {\rm Inv}(w)=\emptyset$,
it follows from $z'\leq_L w$ that $\{(j,i)\mid 1\leq j<i\}\cap {\rm Inv}(z')=\emptyset$, and hence
$k\in[i,n-1]$. If $k\in[i+1,n-1]$, then
$$
(k,n)\in{\rm Nsp}(w)\setminus {\rm Nsp}(z')\quad {\rm and}\quad
\{(j,k)\mid j\in[i]\}\subseteq{\rm Inv}(w)\setminus {\rm Inv}(z').
$$
It follows that $\ell(z')\leq \ell(w)-i-1$, contradicting $z'\in B$.
So we must have $k=i$, that is, $z'_i=n$, and thus $z_i=n$ by the definition of $z$.
Consequently, $z\in A\setminus\{v\}$. It is now obvious that $\psi(z)=z'$, showing that $\psi$ is a surjection and hence a bijection.
This completes the proof.
\end{proof}

In order to get the critical Lemma \ref{notranksymmetric(n-1)}, we give the following technical result.

\begin{lemma}\label{lem:iundernn-1}
Let $w=w_1w_2\cdots w_n\in \mathfrak{B}_n$, where $w_n=n-1$, $w_i=\underline{n}$,
$\{w_1, \ldots, w_{i-1}\}=\{\underline{n-i}, \ldots,\underline{n-2}\}$, and $\{w_{i+1}, \ldots, w_{n-1}\}=\{\underline{1}, \ldots, \underline{n-i-1}\}$ for some $i\in [n-2]$.
Let $u\in [e,w]_L$.
\begin{enumerate}
\item\label{item:unun-1eior1} We have either $u_n=n-1$ and $u_i=\underline{n}$, or $u_n=n$ and $u\leq_L s_{n-1}w$.
\item\label{item:unun-1uw-k2} Suppose $u_n=n-1$ and $u_i=\underline{n}$. If $k,r$ are positive integers such that $\ell(u)=\ell(w)-k$ and $r\leq {\rm min}\{i,n-k-1\}$, then
$\{\underline{n-r}, \ldots, \underline{n-3},\underline{n-2}\}\subseteq \{u_1,u_2, \ldots, u_{i-1}\}$.
\end{enumerate}
\end{lemma}
\begin{proof}
\eqref{item:unun-1eior1}
Since $u\leq_Lw$, it follows from ${\rm Neg}(u)\subseteq{\rm Neg}(w)=[n-1]$ that $u_n>0$, and from
\begin{align}\label{eq:jninvwempty}
{\rm Inv}(u)\cap\{(j,n)\mid j\in[n-1]\}\subseteq{\rm Inv}(w)\cap\{(j,n)\mid j\in[n-1]\}=\emptyset
\end{align}
we see that
$u^{-1}_p<0$ for all $p\in[u_{n}+1,n]$, if $p$ exists.
Note that
\begin{align}\label{eq:jnnspsubeq}
{\rm Nsp}(u)\cap \{(j,n)\mid j\in[n-1]\}\subseteq {\rm Nsp}(w)\cap\{(j,n)\mid j\in[n-1]\}=\{(i,n)\}.
\end{align}
So $u_{n}\geq n-1$, that is, $u_n\in\{n-1,n\}$.

If $u_n=n-1$, then, by Eq.\,\eqref{eq:jninvwempty}, there exists $j\in[n-1]$ such that $u_j=\underline{n}$,
which together with Eq.\,\eqref{eq:jnnspsubeq} implies that $j=i$, that is, $u_{i}=\underline{n}$.
If $u_n=n$, then
${\rm Inv}(u)\subseteq {\rm Inv}(w)={\rm Inv}(s_{n-1}w)$, ${\rm Neg}(u)\subseteq {\rm Neg}(w)={\rm Neg}(s_{n-1}w)$, and ${\rm Nsp}(u)\subseteq {\rm Nsp}(w)\backslash\{(i,n)\}={\rm Nsp}(s_{n-1}w)$,
according to which we obtain that $u\leq_L s_{n-1}w$.

\eqref{item:unun-1uw-k2}
Let $M=\{\underline{n-r}, \ldots, \underline{n-3},\underline{n-2}\}$.
From $u\in[e,w]_L$ we see that $\#{\rm Inv}(w)-\#{\rm Inv}(u)\geq0$.
Since $u_i=w_i=\underline{n}$ and $u_n=w_n=n-1$, it follows from Eq.\,\eqref{eq:eqell(w)} that
\begin{align*}
k=\ell(w)-\ell(u)\geq-\sum_{\{j\in[n]\mid w_j<0\}}w_j+\sum_{\{j\in[n]\mid u_j<0\}}u_j=-\sum_{\{j\in[n-1]\setminus\{i\}\mid w_j<0\}}w_j+\sum_{\{j\in[n-1]\setminus\{i\}\mid u_j<0\}}u_j.
\end{align*}
Note that
$$\{u_j\mid u_j<0,j\in[n-1]\backslash\{i\}\}\subseteq[\underline{n-2},\underline{1}]=\{w_j\mid w_j<0,j\in[n-1]\backslash\{i\}\}.$$
So the largest positive integer among $u_j$ where $j\in[n-1]\setminus\{i\}$ must be less than $k+1$.
But from  $r\leq {\rm min}\{i,n-k-1\}$, we have $k+1\leq n-r$ so that
$$M\subseteq\{u_j\mid u_j<0,j\in[n-1]\backslash\{i\}\}\subseteq \{u_1,\ldots, u_{i-1},u_{i+1},\ldots, u_{n-1}\}.$$
Thus, it suffices to show that $M\cap \{u_{i+1},\ldots, u_{n-1}\}=\emptyset$.
Assume that there exists $q\in[i+1,n-1]$ such that $u_q\in M$. Since $r\leq i$, we have $|u_q|\not\in [n-i-1]$.
But the cardinality of $\{|u_{i-1}|,|u_{i+1}|,\ldots,|u_{n-1}|\}$ is $n-i-1$, so $[n-i-1]\nsubseteqq \{|u_{i-1}|,|u_{i+1}|,\ldots,|u_{n-1}|\}$.
Since $\{|u_1|,\ldots,|u_{i-1}|,|u_{i+1}|,\ldots,|u_{n-1}|\}=[n-2]$, we see that  $[n-i-1]\cap \{|u_{1}|,|u_{2}|,\ldots,|u_{i-1}|\}\neq \emptyset$,
so there exists $q'\in[i-1]$ such that $|u_{q'}|\in [n-i-1]$, and thus $(q',q)\in {\rm Inv}(u)$.
But $(q',q)\not\in {\rm Inv}(w)$, contradicting $x\leq_Lw$.
Thus, $M\subseteq \{u_1,u_2, \ldots, u_{i-1}\}$, as claimed.
\end{proof}

\begin{lemma}\label{notranksymmetric(n-1)}
Let $w=w_1w_2\cdots w_n\in \mathfrak{B}_n$, where $w_n=n-1$ and {$w_i=\underline{n}$ for} some $i\in [n-2]$, and let $f(q)=\Lambda_w^L(q)$. If $w$ and $w^{-1}$ are both minimal non-separable, then
\begin{align*}
[q^d]f=\begin{cases}
[q^{\ell(w)-d}]f, &{\rm if}\  0\leq d\leq i-1,\\
[q^{\ell(w)-d}]f-1, &{\rm if}\ d=i.
\end{cases}
\end{align*}
In particular, $\Lambda^{L}_{w}$  is not rank-symmetric.
\end{lemma}
\begin{proof}
Let $w'=w_1\cdots w_{i-1}w_{i+1}\cdots w_{n-1}$,
and let $h_1(q)=[e,s_{n-1}w]_L(q)$ and $h_2(q)=[e,w']_L(q)$ be the rank generating functions.
Since $w$ is minimal nonseparable and
$$s_{n-1}w=w_1\cdots w_{i-1}(\underline{n-1})w_{i+1}\cdots w_{n-1}n=w_{J_{n-1}},$$
where $J_{n-1}=\Delta\backslash\{\alpha_{n-1}\}$, we see that $s_{n-1}w$ is separable.
But $w'$ is a pattern in $s_{n-1}w$, so $w'$ is also separable by Lemma \ref{sepBnavoidpattern}.
Then, by Lemma \ref{rank-symmetric}, both $h_1(q)$ and $h_2(q)$ are symmetric polynomials.

\begin{figure}[t]
\begin{center}
\begin{tikzpicture}\label{Figure 3.}
\node (underline{2}underline{3}underline{5}underline{1}4) at (0,2.4) {{\tiny$\underline{2}\thinspace\underline{3}\thinspace\underline{5}\thinspace\underline{1}\thinspace 4$}};

\node (underline{2}underline{3}underline{5}14) at (-1.5,1.2) {{\tiny$\underline{2}\thinspace\underline{3}\thinspace\underline{5}\thinspace1\thinspace4$}};
\node (underline{3}underline{2}underline{5}underline{1}4) at (0,1.2) {{\tiny$\underline{3}\thinspace\underline{2}\thinspace\underline{5}\thinspace\underline{1}\thinspace4$}};
\node (underline{2}underline{3}underline{4}underline{1}5) at (1.5,1.2) {{\tiny$\underline{2}\thinspace\underline{3}\thinspace\underline{4}\thinspace\underline{1}\thinspace 5$}};

\node (underline{1}underline{3}underline{5}24) at (-3,0) {{\tiny$\underline{1}\thinspace\underline{3}\thinspace\underline{5}\thinspace2\thinspace 4$}};
\node (underline{3}underline{2}underline{5}14) at (-1.5,0) {{\tiny$\underline{3}\thinspace\underline{2}\thinspace\underline{5}\thinspace1\thinspace4$}};
\node (underline{2}underline{3}underline{4}15) at (0,0) {{\tiny$\underline{2}\thinspace\underline{3}\thinspace\underline{4}\thinspace1\thinspace 5$}};
\node (underline{3}underline{2}underline{4}underline{1}5) at (1.5,0) {{\tiny$\underline{3}\thinspace\underline{2}\thinspace\underline{4}\thinspace\underline{1}\thinspace 5$}};
\node (underline{2}underline{4}underline{3}underline{1}5) at (3,0) {{\tiny$\underline{2}\thinspace\underline{4}\thinspace\underline{3}\thinspace\underline{1}\thinspace 5$}};

\node (1underline{3}underline{5}24) at (-5.25,-1.2) {{\tiny$1\thinspace\underline{3}\thinspace\underline{5}\thinspace2\thinspace 4$}};
\node (underline{1}underline{2}underline{5}34) at (-3.75,-1.2) {{\tiny$\underline{1}\thinspace\underline{2}\thinspace\underline{5}\thinspace3\thinspace 4$}};
\node (underline{3}underline{1}underline{5}24) at (-2.25,-1.2) {{\tiny$\underline{3}\thinspace\underline{1}\thinspace\underline{5}\thinspace2\thinspace4$}};
\node (underline{1}underline{3}underline{4}25) at (-0.75,-1.2) {{\tiny$\underline{1}\thinspace\underline{3}\thinspace\underline{4}\thinspace2\thinspace 5$}};
\node (underline{3}underline{2}underline{4}15) at (0.75,-1.2) {{\tiny$\underline{3}\thinspace\underline{2}\thinspace\underline{4}\thinspace1\thinspace 5$}};
\node (underline{2}underline{4}underline{3}15) at (2.25,-1.2) {{\tiny$\underline{2}\thinspace\underline{4}\thinspace\underline{3}\thinspace1\thinspace 5$}};
\node (underline{4}underline{2}underline{3}underline{1}5) at (3.75,-1.2) {{\tiny$\underline{4}\thinspace\underline{2}\thinspace\underline{3}\thinspace\underline{1}\thinspace 5$}};
\node (underline{3}underline{4}underline{2}underline{1}5) at (5.25,-1.2) {{\tiny$\underline{3}\thinspace\underline{4}\thinspace\underline{2}\thinspace\underline{1}\thinspace 5$}};

\node (1underline{3}underline{4}25) at (-2,-2.4) {{\tiny$1\thinspace\underline{3}\thinspace\underline{4}\thinspace2\thinspace 5$}};
\node (underline{1}underline{2}underline{4}35) at (-0.5,-2.4) {{\tiny$\underline{1}\thinspace\underline{2}\thinspace\underline{4}\thinspace3\thinspace 5$}};
\node (underline{3}underline{1}underline{4}25) at (1,-2.4) {{\tiny$\underline{3}\thinspace\underline{1}\thinspace\underline{4}\thinspace2\thinspace 5$}};
\node (underline{1}underline{4}underline{3}25) at (2.5,-2.4) {{\tiny$\underline{1}\thinspace\underline{4}\thinspace\underline{3}\thinspace2\thinspace 5$}};
\node (underline{4}underline{2}underline{3}15) at (4,-2.4) {{\tiny$\underline{4}\thinspace\underline{2}\thinspace\underline{3}\thinspace1\thinspace 5$}};
\node (underline{3}underline{4}underline{2}15) at (5.5,-2.4) {{\tiny$\underline{3}\thinspace\underline{4}\thinspace\underline{2}\thinspace1\thinspace 5$}};
\node (underline{4}underline{3}underline{2}underline{1}5) at (7,-2.4) {{\tiny$\underline{4}\thinspace\underline{3}\thinspace\underline{2}\thinspace\underline{1}\thinspace 5$}};
\draw[dashed,line width=0.5pt] (0.9,1.7)--(2.1,1.7)--(2.1,0.7)--(0.91,0.7)--(0.9,1.7);
\node (F0) at (3.1,1.2) {{\tiny$F_0$}};
\draw[solid,line width=0.8pt] (2.1,1.2)--(F0);

\draw[dashed,line width=0.5pt] (-0.6,0.5)--(3.6,0.5)--(3.6,-0.5)--(-0.6,-0.5)--(-0.6,0.5);
\node (F1) at (4.6,0) {{\tiny$F_1$}};
\draw[solid,line width=0.8pt] (3.6,0)--(F1);

\draw[dashed,line width=0.5pt] (-1.35,-0.7)--(5.85,-0.7)--(5.85,-1.7)--(-1.35,-1.7)--(-1.35,-0.7);
\node (F2) at (6.85,-1.2) {{\tiny$F_2$}};
\draw[solid,line width=0.8pt] (5.85,-1.2)--(F2);

\draw[dashed,line width=0.5pt] (-2.6,-1.9)--(7.6,-1.9)--(7.6,-2.9)--(-2.6,-2.9)--(-2.6,-1.9);
\node (F3) at (8.6,-2.4) {{\tiny$F_3$}};
\draw[solid,line width=0.8pt] (7.6,-2.4)--(F3);

\draw[dashed,line width=0.5pt] (-0.6,2.8)--(-0.6,2)--(0.6,2)--(0.6,2.8)--(-0.6,2.8);
\node (E0) at (-1.6,2.4) {{\tiny$E_0$}};
\draw[solid,line width=0.8pt] (-0.6,2.4)--(E0);

\draw[dashed,line width=0.5pt] (-2.1,1.7)--(-2.1,0.7)--(0.6,0.7)--(0.6,1.7)--(-2.1,1.7);
\node (E1) at (-3.1,1.2) {{\tiny$E_1$}};
\draw[solid,line width=0.8pt] (-2.1,1.2)--(E1);

\draw[dashed,line width=0.5pt] (-3.6,0.5)--(-3.6,-0.5)--(-0.9,-0.5)--(-0.9,0.5)--(-3.6,0.5);
\node (E2) at (-4.6,0) {{\tiny$E_2$}};
\draw[solid,line width=0.8pt] (-3.6,0)--(E2);

\draw[dashed,line width=0.5pt] (-5.85,-0.7)--(-5.85,-1.7)--(-1.65,-1.7)--(-1.65,-0.7)--(-5.85,-0.7);
\node (E3) at (-6.85,-1.2) {{\tiny$E_3$}};
\draw[solid,line width=0.8pt] (-5.85,-1.2)--(E3);

\draw[solid,line width=0.8pt] (underline{2}underline{3}underline{5}underline{1}4)--(underline{2}underline{3}underline{4}underline{1}5)--(underline{2}underline{4}underline{3}underline{1}5)- -(underline{3}underline{4}underline{2}underline{1}5)
(underline{3}underline{2}underline{5}underline{1}4)--(underline{3}underline{2}underline{4}underline{1}5)--(underline{4}underline{2}underline{3}underline{1}5)- -(underline{4}underline{3}underline{2}underline{1}5)
(underline{2}underline{3}underline{5}14)--(underline{2}underline{3}underline{4}15)--(underline{2}underline{4}underline{3}15)- -(underline{3}underline{4}underline{2}15)
(underline{3}underline{2}underline{5}14)--(underline{3}underline{2}underline{4}15)--(underline{4}underline{2}underline{3}15)
(underline{1}underline{3}underline{5}24)--(underline{1}underline{3}underline{4}25)--(underline{1}underline{4}underline{3}25)
(1underline{3}underline{5}24)--(1underline{3}underline{4}25)
(underline{1}underline{2}underline{5}34)--(underline{1}underline{2}underline{4}35)
(underline{3}underline{1}underline{5}24)--(underline{3}underline{1}underline{4}25) ;
\draw[solid,line width=0.8pt] (underline{2}underline{3}underline{5}underline{1}4)--(underline{2}underline{3}underline{5}14)--(underline{1}underline{3}underline{5}24)- -(1underline{3}underline{5}24)
(underline{2}underline{3}underline{5}underline{1}4)--(underline{3}underline{2}underline{5}underline{1}4)--(underline{3}underline{2}underline{5}14)- -(underline{3}underline{1}underline{5}24)
(underline{2}underline{3}underline{4}underline{1}5)--(underline{2}underline{3}underline{4}15)--(underline{1}underline{3}underline{4}25)- -(1underline{3}underline{4}25)  (underline{2}underline{3}underline{4}underline{1}5)--(underline{3}underline{2}underline{4}underline{1}5)--(underline{3}underline{2}underline{4}15)- -(underline{3}underline{1}underline{4}25)
(underline{2}underline{4}underline{3}underline{1}5)--(underline{2}underline{4}underline{3}15)--(underline{1}underline{4}underline{3}25)
(underline{3}underline{4}underline{2}15)--(underline{3}underline{4}underline{2}underline{1}5)--(underline{4}underline{3}underline{2}underline{1}5)
(underline{4}underline{2}underline{3}underline{1}5)--(underline{4}underline{2}underline{3}15)
(underline{1}underline{3}underline{5}24)--(underline{1}underline{2}underline{5}34)
(underline{2}underline{3}underline{5}14)--(underline{3}underline{2}underline{5}14)
(underline{2}underline{3}underline{4}15)--(underline{3}underline{2}underline{4}15)
(underline{1}underline{3}underline{4}25)--(underline{1}underline{2}underline{4}35) ;

\node(Figure 2.)   at (0,-4) {Figure $2.$ Part of the left weak order of $[e,\underline{2}\thinspace\underline{3}\thinspace\underline{5}\thinspace\underline{1}\thinspace 4]_L$.};
\end{tikzpicture}
\end{center}
\end{figure}

We next divide our proof into three steps. First, we show that for any $d\leq i$,
\begin{align}\label{eq:wn-1qdfh1}
	[q^d]f=[q^{\ell(s_{n-1}w)-d}]h_1
\end{align}
and 	
\begin{align}\label{eq:w-dqdfh1+h2}
[q^{\ell(w)-d}]f=[q^{\ell(s_{n-1}w)-d+1}]h_1+[q^{\ell(w')-d}]h_2.
\end{align}
To see this, let $u\in [e,w]_L$, and let
\begin{align*}
E_d&=\{x\in [e,w]_L\mid  x_n=n-1,\ell(x)=\ell(w)-d\},\\
F_d&=\{x\in [e,s_{n-1}w]_L\mid \ell(x)=\ell(s_{n-1}w)-d\},
\end{align*}
See Figure $2$ for an example in $\mathfrak{B}_5$, where $w=\underline{2}\thinspace\underline{3}\thinspace\underline{5}\thinspace\underline{1}\thinspace 4$.
It follows from Lemma \ref{lem:wjk>0<0uponwn+-} that
\begin{align}\label{eq:w1wi-1i+1n-1}
\{w_1, \ldots, w_{i-1}\}=\{\underline{n-i}, \ldots,\underline{n-2}\}\quad \text{and}\quad \{w_{i+1}, \ldots, w_{n-1}\}=\{\underline{1}, \ldots, \underline{n-i-1}\}.
\end{align}
According to Lemma \ref{lem:iundernn-1}\eqref{item:unun-1eior1}, we have $u_n\in\{n-1,n\}$ so that
\begin{align*}
[q^{\ell(w)-d}]f=\#E_d+\#\{x\in[e,w]_L\mid x_n=n,\ell(x)= \ell(w)-d\}.
\end{align*}
Again by Lemma \ref{lem:iundernn-1}\eqref{item:unun-1eior1}, if $u_n=n$, then $u\leq_L s_{n-1}w$, it follows from $s_{n-1}w\leq_L w$ and $\ell(w)=\ell(s_{n-1}w)+1$ that
$$
\{x\in[e,w]_L\mid x_n=n,\ell(x)= \ell(w)-d\}=F_{d-1}.
$$
Hence, $$\#\{x\in[e,w]_L\mid x_n=n,\ell(x)= \ell(w)-d\}=[q^{\ell(s_{n-1}w)-d+1}]h_1.$$
If $u_n=n-1$, then, by Lemma \ref{lem:iundernn-1}\eqref{item:unun-1eior1}, $u_i=\underline{n}$,
and by Lemma \ref{leftweakorder}, $u'\leq_L w'$ in $\mathfrak{B}_{n-2}$, where $u'=u_1\cdots u_{i-1}u_{i+1}\cdots u_{n-1}$.
In other words,
$u\in[e,w]_L$ with $u_{n}=n-1$ is equivalent to
$u'\leq_L w'$ in $\mathfrak{B}_{n-2}$.
Since by Eq.\,\eqref{eq:eqell(w)} $\ell(w)-\ell(u)=\ell(w')-\ell(u')$ , we see that the map
\begin{align*}%\label{eq:xn=n-1bijecvarphi}
\varphi:E_d\rightarrow\{x\in [e,w']_L\mid \ell(x)= \ell(w')-d\},\quad x\mapsto x'=x_1\cdots x_{i-1}x_{i+1}\cdots x_{n-1}
\end{align*}
is a bijection.
So
\begin{align}\label{eq:edqelw-dh2}
\#E_d=[q^{\ell(w')-d}]h_2,
\end{align}
and hence Eq.\,\eqref{eq:w-dqdfh1+h2} follows.

Assume that $u\in[e,w]_L$ with $\ell(u)=d$. Using Lemma \ref{lem:iundernn-1}\eqref{item:unun-1eior1},
if $u_{n}=n-1$ and $u_i=\underline{n}$, then $\ell(u)\geq n>i\geq d$ by Eq.\,\eqref{eq:eqell(w)}, {contradicting}  $\ell(u)=d$.
Thus, we must have $u_n=n$ and $u\leq s_{n-1}w$. Combining this with $s_{n-1}w\leq_Lw$ gives that
$$
\{x\in[e,w]_L\mid \ell(x)= d\}=\{x\in[e,s_{n-1}w]\mid \ell(x)=d\},
$$
that is, $[q^d]f=[q^d]h_1$ for all $d\leq i$. But $h_1(q)$ is symmetric, so Eq.\,\eqref{eq:wn-1qdfh1} holds.

The next thing to do in the proof is to show that
\begin{align}\label{eq:asn-1wh12}
[q^{\ell(s_{n-1}w)-d}]h_1=\begin{cases}
\sum\limits^d_{j=0}[q^{\ell(w')-j}]h_2,&{\rm if}\ 0\leq d<i,\\
\sum\limits^d_{j=0}[q^{\ell(w')-j}]h_2-1,&{\rm if}\  d=i.
\end{cases}
\end{align}
To this end, for each $d\leq i$, let
\begin{align*}
G_d=\begin{cases}
\biguplus_{j=0}^dE_j, &{\rm if}\ 0\leq d<i,\\
\left(\biguplus_{j=0}^iE_j\right)\backslash\{w\}, &{\rm if}\ d=i.
\end{cases}
\end{align*}
It is clear that $\{w\}= G_0\subseteq G_1\subseteq \cdots\subseteq G_{i-1}\subseteq G_{i}\cup\{w\}$.
For each $d\leq i$, define a map
\begin{align*}
\xi_d: G_d \rightarrow F_d,\quad x \mapsto  s_{n-p_d}\cdots s_{n-2}s_{n-1}x,
\end{align*}
where $p_d=\ell(x)-\ell(s_{n-1}w)+d$. Note that $p_d$ is indeed determined by $x$ and $d$.
To prove Eq.\,\eqref{eq:asn-1wh12} it suffices, in view of Eq.\,\eqref{eq:edqelw-dh2}, to show that $\xi_d$ are bijections for each $d\leq i$.
This is trivial if $d=0$ since $G_0=\{w\}$ and $\xi_0(w)=s_{n-1}w$ is the unique element of $F_0$. We next suppose $d\in[i]$.

To show that $\xi_d$ is well-defined, we first prove $\xi_d(w)\in F_d$ for $d<i$.
Assume that $0<d<i$. Then
$$p_d=\ell(w)-\ell(s_{n-1}w)+d= d+1\in[i],$$
and hence $\xi_d(w)=s_{n-d-1}\cdots s_{n-2}s_{n-1}w$.
It follows from Eq.\,\eqref{eq:w1wi-1i+1n-1} that
$\ell(\xi_d(w))=\ell(s_{n-1}w)-d$, that is, $\xi_d(w)\in F_d$.
It remains to show that $\xi_d(x)\in F_d$ for any $x\in G_d$ with $x\neq w$, where $d\leq i$.
To do so we take $x\in G_d\backslash\{w\}$ and set $k=\ell(w)-\ell(x)$.
Then $k\in[d]$, $x\in E_k$ and $p_d=d-k+1$, so that
$\xi_d(x)=s_{n-d+k-1}\cdots s_{n-2}s_{n-1}x$. By Lemma \ref{lem:iundernn-1} and Eq.\,\eqref{eq:w1wi-1i+1n-1}, we have $x_n=n-1$ and $x_i=\underline{n}$. Combining this with
Lemma \ref{leftweakorder} we get that
$s_{n-1}x\leqslant_Ls_{n-1}w$. If $k=d$, then $x\in E_d$ and $\xi_d(x)=s_{n-1}x$, which together with
\begin{align*}
\ell(s_{n-1}x)=\ell(x)-1=\ell(w)-d-1=\ell(s_{n-1}w)-d
\end{align*}
yields that $s_{n-1}x\in F_d$, that is, $\xi_d(x)\in F_d$. If $k\in [d-1]$, then $2\leq d-k+1\leq d\leq i\leq n-2$, so $d-k+1 \leq {\rm min}\{i,n-k-1\}$.
Applying Lemma \ref{lem:iundernn-1}\eqref{item:unun-1uw-k2}, we obtain that $\{\underline{n-d+k-1}, \ldots, \underline{n-3},\underline{n-2}\}\subseteq \{x_1,x_2, \ldots, x_{i-1}\}$.
So
$$
\ell(s_{n-j-1}s_{n-j}\cdots s_{n-1}x)=\ell(s_{n-j}\cdots s_{n-1}x)-1
$$
for all $j\in[d-k]$. Consequently, $\xi_{d}(x)\leq_L s_{n-1}w$ and
$$\ell(\xi_{d}(x))=\ell(s_{n-1}x)-d+k=\ell(x)-d+k-1=\ell(w)-d-1=\ell(s_{n-1}w)-d.$$
So $\xi_{d}(x)\in F_d$, and we conclude that $\xi_d$ is well defined.

To show that $\xi_d$ is injective, we take two elements $x$ and $y$ in $G_d$ and assume that $\xi_d(x)=\xi_d(y)$.
Then there exist $r,t\in [d]$ such that $x\in E_r$ and $y\in E_t$.
Without loss of generality we may assume that $r\leq t$.
By the definition of $\xi_d$, we see that
\begin{align}\label{eq:xidx=xidy}
s_{n-d+r-1}\cdots s_{n-1}x=\xi_d(x)= \xi_d(y)=s_{n-d+t-1}\cdots s_{n-1}y.
\end{align}
It follows from Lemma \ref{lem:iundernn-1}\eqref{item:unun-1uw-k2} that
$\{\underline{n-d+r-1}, \ldots, \underline{n-3},\underline{n-2}\}\subseteq \{x_1,x_2, \ldots, x_{i-1}\}$,
from which we get
${\rm Neg}(x)={\rm Neg}(\xi_d(x))$, ${\rm Nsp}(x)={\rm Nsp}(\xi_d(x))\uplus\{(i,n)\}$ and
$${\rm Inv}(x)={\rm Inv}(\xi_d(x))\uplus\{(j,i)\mid j\in[i-1],x_j\in\{\underline{n-d+r-1}, \ldots, \underline{n-3},\underline{n-2}\}\},$$
and similarly for $y$. Since $r\leq t$, it follows from Theorem \ref{leftweakorder} that $y\leq_L x$.
From $x\in E_r$ and $y\in E_t$ we see that $\ell(x)-\ell(y)=t-r$.
So there exist simple reflections $s'_1, \ldots, s'_{t-r}$ such that $s'_1\cdots s'_{t-r}y=x$. Substituting into Eq.\,\eqref{eq:xidx=xidy} gives
$$
s_{n-d+r-1}\cdots s_{n-1}=s_{n-d+t-1}\cdots s_{n-1}s'_{t-r}\cdots s'_{1},
$$
which are reduced.
But all the braid-moves in $\mathfrak{B}_n$ are $s_0s_1s_0s_1=s_1s_0s_1s_0$,
 $s_is_{i+1}s_i=s_{i+1}s_{i}s_{i+1}$ and $s_is_j=s_js_i$, where $1\leq i<j\leq n-1$ with $j\neq i+1$.
Hence, the above reduced expression $s_{n-d+r-1}\cdots s_{n-1}$ cannot be changed by braid-moves. However, every two reduced expressions for an element can be connected via a sequence of braid-moves (see, for example, \cite[Theorem 3.3.1]{BB05}),
then $s_{n-d+r-1}\cdots s_{n-1}$ has a unique reduced expression. So we must have $r=t$.
Thus, $x=y$ follows from Eq.\,\eqref{eq:xidx=xidy}, showing that $\xi_d$ is injective.

To show $\xi_d$ is surjective, we take any $z\in F_d$.
{Since} $\ell(s_{n-1}w)-\ell(z)=d\leq n-2$, it follows from Eq.\,\eqref{eq:eqell(w)} that $z^{-1}_{n-1}<0$.
By Lemma \ref{lem:iundernn-1}, we have $z_n=n$, and hence $|z_i|\leq n-1$.
Let $x=s_{n-1}s_{n-2}\cdots s_{|z_i|}z$.
It suffices to show that $x\in G_d$ and $\xi_d(x)=z$.
If $z_i=\underline{n-1}$, then $x=s_{n-1}z$,
which is an element of $[e,w]_L$ with $x_n=n-1$ and $\ell(x)=\ell(z)+1=\ell(w)-d$, so $x\in E_d\subseteq G_d$ and $\xi_d(x)=z$.
If $z_i\geq\underline{n-2}$, we claim that $z_i<0$ and
\begin{align}\label{eq:n-1subseteqz-i}
\{\underline{n-1}, \underline{n-2}, \ldots,z_i-1\}\subseteq \{z_1,z_2, \ldots, z_{i-1}\}.
\end{align}
Since $z\leq_Ls_{n-1}w$, it follows from Eq.\,\eqref{eq:w1wi-1i+1n-1} that
\begin{align}\label{eq:invsetsn-1w-z}
\{(i,j)\mid j\in[i+1,n-1]\}\cap {\rm Inv}(z)\subseteq \{(i,j)\mid j\in[i+1,n-1]\}\cap {\rm Inv}(s_{n-1}w)=\emptyset,
\end{align}
If $z_i>0$, then $z_j>z_i>0$ for all $j\in[i+1,n-1]$.
Hence
\begin{align*}
{\rm Inv}(s_{n-1}w)\setminus{\rm Inv}(z)&\supseteq\{(j,i)\mid j\in[i-1], z_j<0\},\\
{\rm Neg}(s_{n-1}w)\setminus{\rm Neg}(z)&=[i,n-1]\cup\{j\in[i-1]\mid z_j>0\},
\end{align*}
which together with Eq.\,\eqref{eq:eqell(w)} yields that
$$
\ell(s_{n-1}w)-\ell(z)\geq \#{\rm Inv}(s_{n-1}w)-\#{\rm Inv}(z)+\#{\rm Neg}(s_{n-1}w)-\#{\rm Neg}(z)\geq n-1>d,
$$
contradicting $z\in F_d$. So $z_i<0$.
We show Eq.\,\eqref{eq:n-1subseteqz-i} holds by induction.
Since $z^{-1}_{n-1}<0$, it follows from Eq.\,\eqref{eq:invsetsn-1w-z} that
$\underline{n-1}\in \{z_1,z_2,\ldots,z_{i-1}\}$.
Suppose $t\in [\underline{n-1},z_i-1]$ and $[\underline{n-1},t-1]\subseteq\{z_1,z_2,\ldots,z_{i-1}\}$.
Comparing the inversion sets of $z$ and $s_{n-1}w$, it follows from Eq.\,\eqref{eq:w1wi-1i+1n-1} that
\begin{align*}
{\rm Inv}(s_{n-1}w)-{\rm Inv}(z)\geq \#\{(j,i)\mid j\in[i-1], z_j\in [\underline{n-1},t-1]\}=n+t-1.
\end{align*}
Let $j=|z^{-1}_t|$, that is, $z_j\in\{t,\underline{t}\}$.
If $z_j=\underline{t}$, then since $t$ is a negative entry of $s_{n-1}w$, from Eq.\,\eqref{eq:w1wi-1i+1n-1} we see that
\begin{align*}
\ell(s_{n-1}w)-\ell(z)\geq \#{\rm Inv}(s_{n-1}w)-\#{\rm Inv}(z)-t\geq n-1,
\end{align*}
contradicting $z\in F_d$. This means that $z_j=t$, so that $z_j<z_i$.
By Eq.\,\eqref{eq:invsetsn-1w-z}, we have $j<i$, that is, $t\in\{z_1,z_2,\ldots,z_{i-1}\}$.
So Eq.\,\eqref{eq:n-1subseteqz-i} follows by induction.
Comparing the cardinalities of the two sets in Eq.\,\eqref{eq:n-1subseteqz-i} yields that $n-1+z_i\leq i-1$, that is, $z_i\leq \underline{n-i}$.
Thus, $\underline{n-2}\leq z_i\leq \underline{n-i}$.
It follows Eq.\,\eqref{eq:n-1subseteqz-i} that
$x_{n}=n-1$ and
\begin{align*}
\ell(x)=\ell(z)+n+z_i=\ell(s_{n-1}w)-d+n+z_i=\ell(w)-1-d+n+z_i,
\end{align*}
so $\ell(w)-\ell(x)=1+d-n-z_i\leq d-1$.
Moreover, we have ${\rm Neg}(x)={\rm Neg}(z)\subseteq{\rm Neg}(w)$,
${\rm Nsp}(x)={\rm Nsp}(z)\cup\{(i,n)\}\subseteq{\rm Nsp}(w)$, and
$$
{\rm Inv}(x)={\rm Inv}(z)\cup\{(j,i)\mid j\in[i-1],z_j\leq z_i-1\}\subseteq{\rm Inv}(z)\cup\{(j,i)\mid j\in[i-1]\}\subseteq{\rm Nsp}(w).
$$
Hence, $x\leq_Lw$ so that $x\in G_d$. It is now routine to verify that $\xi_d(x)=z$, showing that $\xi_d$ is surjective and hence is bijective.
So Eq.\,\eqref{eq:asn-1wh12} is true.

Finally, we have to show the desired identity holds.
By Eq.\,\eqref{eq:wn-1qdfh1}, we have $[q^d]f=[q^{\ell(s_{n-1}w)-d}]h_1$ for all $d\leq i$.
Combining this with Eqs.\,\eqref{eq:w-dqdfh1+h2} and \eqref{eq:asn-1wh12} yields that for $d<i$, we have
\begin{align*}
[q^d]f=\sum\limits^{d-1}_{j=0}[q^{\ell(w')-j}]h_2+[q^{\ell(w')-d}]h_2
=[q^{\ell(s_{n-1}w)-d+1}]h_1+[q^{\ell(w')-d}]h_2
=[q^{\ell(w)-d}]f,
\end{align*}
and for $d=i$, we have
\begin{align*}
[q^d]f&=\sum\limits^{d-1}_{j=0}[q^{\ell(w')-j}]h_2+[q^{\ell(w')-d}]h_2-1\\
&=[q^{\ell(s_{n-1}w)-d+1}]h_1+[q^{\ell(w')-d}]h_2-1\\
&=[q^{\ell(w)-d}]f-1.
\end{align*}
This completes the proof.
\end{proof}

Before proceeding further, let's briefly explain the relationship among the rank generating functions of $[e,u]_L$, $[e,u^{-1}]_L$, $[e,u]_R$, and $[e,u^{-1}]_R$ for any $u\in \mathfrak{B}_n$.
Since $x\mapsto xu^{-1}$ is a bijection from $[e,u]_L$ to $[e,u^{-1}]_L$ such that $\ell(xu^{-1})=\ell(u)-\ell(x)$, we conclude that $[e,u]_L$ and $[e,u^{-1}]_L$ have the same rank-symmetry and rank-unimodality.
Since the left and right weak orders are isomorphic via the map $w\mapsto w^{-1}$, we have $[e,u]_L=[e,u^{-1}]_R^{-1}$, and similarly for $[e,u]_R=[e,u^{-1}]_L^{-1}$.
So the rank-generating functions of the four sets have the same symmetry and unimodality.

{\begin{coro}\label{cor:ewrsym}
Let $w=w_1w_2\cdots w_n$ be an element of $\mathfrak{B}_n$ such that both $w$ and $w^{-1}$ are minimal non-separable. If $\{|w_n|,|w_{n-1}|\}\neq\{n-1,n\}$, then $[e,w]_R$ is not rank-symmetric.
\end{coro}
\begin{proof}

Since $w$ and $w^{-1}$ are both minimal non-separable, it follows from Lemma \ref{nwn-1nn-1n} that either $|w_n|=n-1$ or $|w_{n-1}|=n$.
But $\{|w_n|,|w_{n-1}|\}\neq\{n-1,n\}$, so if $|w_{n-1}|= n$, then there exists $i\in[n-2]$ such that $|w_{n}|=i$, and hence $|w_{i}^{-1}|=n$.
By Lemma \ref{minonsepBn}\eqref{minonsepBn:item2}, either $w_{i}^{-1}=n$ and $w_n^{-1}=\underline{n-1}$, or $w_{i}^{-1}=\underline{n}$ and $w_n^{-1}=n-1$.
Replacing the $w$ by $w^{-1}$ in Lemmas \ref{notranksymmetric-(n-1)} and \ref{notranksymmetric(n-1)} yields that $[e,w^{-1}]_L$
and hence $[e,w]_R$ is not rank-symmetric.

If $|w_n|=n-1$, then the above argument applies as well, from Lemmas \ref{notranksymmetric-(n-1)} and \ref{notranksymmetric(n-1)} we see that $[e,w]_L$ and hence $[e,w]_R$ is not rank-symmetric.
\end{proof}

We are now in a position to prove  Theorem \ref{maintheorem}.
To begin with, let's recall two methods of producing new splittings from a given one in Weyl groups.
\begin{lemma}[\cite{GG20am}, Propositions 9 and 10]\label{lemproducesplitting}
Let $(X, Y)$ be a splitting of a Weyl group $W$.
\begin{enumerate}
%\item\label{lemproducesplitting:itema} $X$ and $Y$ have unique maximal elements $x_0$ and $y_0$ under left and right weak order, respectively. Furthermore, we have $x_0y_0=w_0$.
\item\label{lemproducesplitting:itemb} For any $J\subseteq\Delta$, the pair $(X\cap W_J,Y\cap W_J)$ is a splitting of $W$.
\item\label{lemproducesplitting:itemc} $X$ and $Y$ have unique maximal elements $x_0$ and $y_0$ under left and right weak order, respectively. Furthermore, we have $x_0y_0=w_0$.
If we define maps $\varphi:W\rightarrow W$ and $\psi:W\rightarrow W$ by
$$\varphi(z)=zx_0^{-1}\quad \text{and}\quad \psi(z)=x_0zw_0,$$
respectively, then $(\varphi(X), \psi(Y))$ is also a splitting of $W$, having left- and right-maximal elements $x_0^{-1}$ and $x_0w_0$, respectively.
\end{enumerate}
\end{lemma}
}
%For the symmetric groups  Gaetz and Gao \cite[Theorem 3]{GG20am} proved the following result.

We are now turning to the proof of Theorem \ref{maintheorem}.

\begin{proof}[Proof of Theorem \ref{maintheorem}]
Let $W$ be the Weyl group of type $\mathfrak{B}_n$, and let $U=[e,u]_R\subseteq W$.
Then $W/U=[e, w_0u^{-1}]_L$.
According to \cite[Theorem 3]{GG20am}, if $u$ is separable, then $(W/U,U)$ is a splitting of $W$. So it suffices to show the converse is true.
Assume that $(W/U,U)$ is a splitting of $W$.
Then the multiplication map $W/U\times U\rightarrow W$ is length-additive and bijective, so we have the factorization $W(q)=(W/U)(q)U(q)$.
By Eq.\,\eqref{eq:W(q)}, $W(q)$ is a product of cyclotomic polynomials which are irreducible over $\mathbb{Q}$, it follows that $U(q)$ is a product of cyclotomic polynomials and hence is symmetric.

Suppose $u$ is non-separable. {
Then, by Lemma \ref{sepBnavoidpattern}, $u^{-1}$, $w_0u$, and $w_0u^{-1}$ are non-separable.
We first prove by induction on $n$ that both $u$ and $u^{-1}$ are minimal non-separable.
The base case $n = 2$ being trivial by Lemma \ref{minonsepBn}.
Suppose $n\geq 3$, and let $J=\Delta\backslash\{\alpha_i\}$ for some $i\in[0,n-1]$. By Lemma \ref{lemproducesplitting}\eqref{lemproducesplitting:itemb}, the pair $((W/U)\cap W_J, U\cap W_J)$ is a splitting of $W_J$.
Since $U=[e,u]_R\subseteq W$, we conclude that
$$U\cap W_J=[e,u]_R\cap W_J=[e,u^{-1}]_L^{-1}\cap W_J=([e,u^{-1}]_L\cap W_J)^{-1}=[e,(u^{-1})_J]_L^{-1}=[e,((u^{-1})_J)^{-1}]_R.$$
It follows from $W/U=[e, w_0u^{-1}]_L$ that
$$(W/U)\cap W_J=[e, w_0u^{-1}]_L\cap W_J=[e, (w_0u^{-1})_J]_L.$$
Since
\begin{align*}
I_{\Phi_J}((w_0u^{-1})_J)=I_{\Phi}(w_0u^{-1})\cap \Phi_J^{+}=\left(\Phi^{+}\setminus I_{\Phi}(u^{-1})\right)\cap \Phi_J^{+}=\Phi_J^{+}\setminus I_{\Phi}(u^{-1})=I_{\Phi_J}(w_0(J)(u^{-1})_{J}),
\end{align*}
we see that $(w_0u^{-1})_J=w_0(J)(u^{-1})_{J}$, and hence
$$(W/U)\cap W_J=[e, w_0(J)(u^{-1})_{J}]_L=W_J/(U\cap W_J).$$
If $i=0$, then $W_J=\mathfrak{S}_n$, thus, by \cite[Theorem 4]{GG20am}, $((u^{-1})_J)^{-1}$ and hence $(u^{-1})_J$ is separable. If $i\in[n-1]$, then
$W_J=\mathfrak{B}_i\times \mathfrak{S}_{n-i}$, so $((u^{-1})_J)^{-1}=u'\times u''$ for some $u'\in\mathfrak{B}_i$ and $u''\in \mathfrak{S}_{n-i}$, from which we get
$U\cap W_J=[e,u']_R\times [e,u'']_R$, so that $(W/U)\cap W_J=\left(\mathfrak{B}_i/[e,u']_R\right)\times \left(\mathfrak{S}_{n-i}/ [e,u'']_R\right)$.
Thus, $(\mathfrak{B}_i/[e,u']_R,[e,u']_R)$ and $(\mathfrak{S}_{n-i}/ [e,u'']_R,[e,u'']_R)$ are splittings of $\mathfrak{B}_{i}$ and $\mathfrak{S}_{n-i}$, respectively.
By the inductive hypothesis, $u'$ is separable, and  by \cite[Theorem 4]{GG20am} again, $u''$ is separable, showing that $(u^{-1})_J$ is separable.
So $(u^{-1})_J$ is separable for each $J=\Delta\backslash\{\alpha_i\}$ with $i\in[0,n-1]$, and hence $u^{-1}$ is minimal non-separable.

Since $(W/U,U)=([e, w_0u^{-1}]_L, [e,u]_R)$ is a splitting of $W$, it follows from Lemma \ref{lemproducesplitting}\eqref{lemproducesplitting:itemc} that $([e, uw_0]_L, [e,w_0u^{-1}w_0]_R)$ is also a splitting of $W$. Then, by \ref{lemproducesplitting}\eqref{lemproducesplitting:itemb},
$$([e, uw_0]_L\cap W_J, [e,w_0u^{-1}w_0]_R\cap W_J)$$ is a splitting of $W_J$ for each $J=\Delta\setminus\{\alpha_i\}$, where $i\in[0,n-1]$.
But for $\mathfrak{B}_n$ we have $uw_0=w_0u$ and $w_0u^{-1}w_0=u^{-1}$, so
$$[e, uw_0]_L\cap W_J=[e, w_0u]_L\cap W_J=[e, (w_0u)_{J}]_L=[e, w_0(J)u_{J}]_L,$$
and
%\begin{align*}
%[e,w_0u^{-1}w_0]_R\cap W_J&=[e,u^{-1}]_R\cap W_J=[e,u]_L^{-1}\cap W_J\\
%&=([e,u]_L\cap W_J)^{-1}\\
%&=[e,u_J]_L^{-1}=[e,u_J^{-1}]_R.
%\end{align*}
$$[e,w_0u^{-1}w_0]_R\cap W_J=[e,u^{-1}]_R\cap W_J=[e,u]_L^{-1}\cap W_J=([e,u]_L\cap W_J)^{-1}=[e,u_J]_L^{-1}=[e,(u_J)^{-1}]_R.$$
Completely analogous to the proof of the minimality of $u^{-1}$, we can show that $u$ is also minimal non-separable.
%It follows from Lemma \ref{ww0w-1separableprop} that $u$ and $u^{-1}$ are also minimal non-separable.

Now we see that $u$ and $u^{-1}$ are both minimal non-separable, so by Lemma \ref{nwn-1nn-1n}, either $|u_n|=n-1$ or $|u_{n-1}|=n$.
If $\{|u_n|,|u_{n-1}|\}=\{n-1,n\}$, then $|u_n|=n-1$ and $|u_{n-1}|=n$. By Lemma \ref{minonsepBn}\eqref{minonsepBn:item2}, $u_{n-1}u_n$ is either $\underline{n}(n-1)$ or $n(\underline{n-1})$.
Let $J=\Delta\backslash\{\alpha_{n-2},\alpha_{n-1}\}$. Since $U=[e,u]_R$, it follows that
\begin{align*}
u^J=12\cdots(n-2)u_{n-1}u_n\in U\cap W/U,
\end{align*}
so that the distinct pairs $(e,u^J)$ and $(u^J,e)$ are both sent to $u^J$ by the multiplication map $W/U \times U\rightarrow W$, contradicting the fact that $(W/U,U)$ is a splitting of $W$.
If $\{|u_n|,|u_{n-1}|\}\neq\{n-1,n\}$, then, by Corollary \ref{cor:ewrsym}, $U=[e,u]_R$ is not rank-symmetric, a contradiction.
Thus, $u$ must be separable, completing the proof.}
\end{proof}

\delete{
Consider the symmetry of the rank-generating function of $\Lambda^{L}_{w}$, where  $w\in\mathfrak{B}_n$ such that both $w$ and $w^{-1}$ are minimal non-separable.
It follows from Lemma \ref{nwn-1nn-1n} that either $|w_n|=n-1$ or $|w_{n-1}|=n$.
With $\{|w_n|,|w_{n-1}|\}\neq\{n-1,n\}$, we can conclude that $\Lambda^{L}_{w}$  is not rank-symmetric.
%We can conclude from above that $\Lambda^{L}_{w}$  is not rank-symmetric if
%$\{|w_n|,|w_{n-1}|\}\neq\{n-1,n\}$.
What's more, if $\{|w_n|,|w_{n-1}|\}=\{n-1,n\}$, then $\Lambda^{L}_{w}$ is rank-symmetric  and rank-unimodal, see Proposition \ref{ranksymmetric}.
Use the two lemmas below to prove this result.}

Let $w\in\mathfrak{B}_n$ such that $w$ and $w^{-1}$ are both minimal non-separable.
It follows from Lemma \ref{nwn-1nn-1n} that either $|w_n|=n-1$ or $|w_{n-1}|=n$.
{According to Corollary \ref{cor:ewrsym}, we see that $\Lambda^{L}_{w}$ is not rank-symmetric if
$\{|w_n|,|w_{n-1}|\}\neq\{n-1,n\}$. However, if $\{|w_n|,|w_{n-1}|\}=\{n-1,n\}$, then $\Lambda^{L}_{w}$ is rank-symmetric  and rank-unimodal,
see Proposition \ref{ranksymmetric}, where only $w$ is required to be minimal non-separable.
We need the following two lemmas to prove this result.}

\begin{lemma}\label{lem:uniqreducedexpre}
Let $w=12\cdots(n-2)\underline{n}(n-1)$. Then $s_{n-1}s_{n-2}\cdots s_1s_0s_1\cdots s_{n-3}s_{n-2}$ is the unique reduced expression of $w$.
\end{lemma}
\begin{proof}
It is straightforward to check that $w=s_{n-1}s_{n-2}\cdots s_1s_0s_1\cdots s_{n-3}s_{n-2}$. It follows from  Eq.\,\eqref{eq:eqell(w)} that $\ell(w)=2n-2$, so this expression is reduced.
Since the above reduced expression of $w$ cannot be changed by braid-moves, by \cite[Theorem 3.3.1]{BB05}, we have $s_{n-1}s_{n-2}\cdots s_1s_0s_1\cdots s_{n-3}s_{n-2}$ is the unique reduced expression of $w$, and the proof follows.
\end{proof}

\begin{lemma}\label{varphilength-additive-bijectionmap}
Let $w=w_1w_2\cdots w_n\in \mathfrak{B}_n$ with $w_{n-1}w_{n}=\underline{n}(n-1)$. Then the map
\begin{align*}
\varphi:\Lambda^{L}_{w^J}\times\Lambda^{L}_{w_J}\rightarrow\Lambda^{L}_w,\quad (x,y)\mapsto xy
\end{align*}
 is a length-additive bijection, where $J=\Delta\backslash\{\alpha_{n-2},\alpha_{n-1}\}$.
\end{lemma}
\begin{proof}
We first show that $\varphi$ is well defined. To this end, let
$x\in\Lambda^{L}_{w^J}$ and $y\in\Lambda^{L}_{w_J}$.
Since $J=\Delta\backslash\{\alpha_{n-2},\alpha_{n-1}\}$, we have
$w_J=w_1w_2\cdots w_{n-2}(n-1)n$ and $w^J=12\cdots(n-2)\underline{n}(n-1).$
Thus, $w=w_Jw^J=w^Jw_J$, and hence $w^J\leq _Lw$ and $w_J\leq_L w$.
Moreover, we can write $y=y_1y_2\cdots y_{n-2}(n-1)n$.
\delete{Moreover,
\begin{align*}
{\rm Neg}(w)={\rm Neg}(w_J)\uplus {\rm Neg}(w^J),\
{\rm Nsp}(w)={\rm Nsp}(w_J)\uplus {\rm Nsp}(w^J),\
{\rm Inv}(w)={\rm Inv}(w_J)\uplus {\rm Inv}(w^J).
\end{align*}}
By Lemma \ref{lem:uniqreducedexpre}, $w^J$ has the unique reduced expression
$w^J=s_{n-1}s_{n-2}\cdots s_1s_0s_1\cdots s_{n-3}s_{n-2},$ so
\begin{align*}
\Lambda^{L}_{w^J}
&=\{e\}\uplus\{s_is_{i+1}\cdots s_{n-3}s_{n-2}\mid i\in[0,n-2]\}\uplus\{s_{i}s_{i-1}\cdots s_1s_0s_1\cdots s_{n-3}s_{n-2}\mid i\in[n-1]\}\\
&=\{e,w^J\}\uplus\{1\cdots \hat{i}\cdots(n-1)\,i\,n\mid i\in[n-2]\}\uplus\{1\cdots\hat{i}\cdots(n-1)\,\underline{i}\,n \mid i\in[n-1]\}.
\end{align*}
If $x=e$, then $xy=y\leq_Lw_J\leq_Lw$.
If $x=w^J$, then
$$xy=\big(1\cdots(n-2)\underline{n}(n-1)\big)\cdot \big(y_1\cdots y_{n-2}(n-1)n\Big)=y_1\cdots y_{n-2}\underline{n}(n-1).$$
Consequently,
\begin{align*}
{\rm Neg}(xy)&={\rm Neg}(y)\uplus {\rm Neg}(w^J)\subseteq{\rm Neg}(w_J)\uplus {\rm Neg}(w^J)={\rm Neg}(w),\\
{\rm Nsp}(xy)&={\rm Nsp}(y)\uplus {\rm Nsp}(w^J)\subseteq{\rm Nsp}(w_J)\uplus {\rm Nsp}(w^J)={\rm Nsp}(w),\\
{\rm Inv}(xy)&={\rm Inv}(y)\uplus {\rm Inv}(w^J)\subseteq{\rm Inv}(w_J)\uplus {\rm Inv}(w^J)={\rm Inv}(w),
\end{align*}
from which we get $xy\leq_L w$.
If $x=1\cdots\hat{i}\cdots(n-1)\,i\,n$ for some $i\in[n-2]$, then
$xy=y'_1\cdots y'_{n-2}\,i\,n$, where
for $j\in[n-2]$, we have
\begin{align*}
y'_j=\begin{cases}
y_j,&{\rm if}\ |y_j|<i,\\
{\rm sgn}(y_j)(|y_j|+1),&{\rm if}\ |y_j|\geq i.
\end{cases}
\end{align*}
Since $y\leq_L w_J\leq_Lw$, there follows
\begin{align*}
{\rm Neg}(xy)&={\rm Neg}(y)\subseteq {\rm Neg}(w),\\
{\rm Nsp}(xy)&\subseteq{\rm Nsp}(y)\cup{\rm Nsp}(w^J)\subseteq {\rm Nsp}(w),\\
{\rm Inv}(xy)&\subseteq{\rm Inv}(y)\cup{\rm Inv}(w^J)\subseteq {\rm Inv}(w),
\end{align*}
from which we conclude that $xy\leq_L w$.
A completely analogous argument shows that $xy\leq_L w$
if $x=1\cdots\hat{i}\cdots(n-1)\,\underline{i}\,n$ for some $i\in[n-1]$.
Therefore, the map $\varphi$  is well defined.

Now we show that $\varphi$ is surjective.  Let $u\in \Lambda^{L}_{w}$. It suffices to show that $u^J\leq_L{w^J}$ and $u_J\leq_L{w_J}$.
Since $u_J\leq_L u\leq_L w$, there follows that $w=v^Jv_Ju_J$, where $v=w(u_J)^{-1}$.
Thus, $u_J\leq_L v_Ju_J=w_J$.
On the other hand, we have $u^J=u_1'u_2'\cdots u_{n-2}'u_{n-1}u_n$, where $u_1'u_2'\cdots u_{n-2}'$ is obtained by rearranging the elements of $\{|u_1|,|u_2|,\ldots,|u_{n-2}|\}$ in increasing order.
Since $u\leq_L w$ and $w_{n}=n-1>0$, we have $u_n>0$ and hence ${\rm Neg}(u^J)\subseteq\{n-1\}={\rm Neg}(w^J)$ and ${\rm Nsp}(u^J)\subseteq{\rm Nsp}(w^J)$.
From ${\rm Inv}(u)\subseteq{\rm Inv}(w)$ we see that $u_{n-1}<u_{n}$, so if ${\rm Inv}(u^J)\nsubseteq{\rm Inv}(w^J)$, then there exists $k\in[n-2]$ such that $u_k'>u_n$.
But $u_k'=|u_{j}|$ for some $j\in[n-2]$, so we have either $(j,n)\in {\rm Nsp}(u)\subseteq{\rm Nsp}(w)$ or $(j,n)\in {\rm Inv}(u)\subseteq{\rm Inv}(w)$, contradicting $w_{n-1}w_{n}=\underline{n}(n-1)$. It follows that ${\rm Inv}(u^J)\subseteq{\rm Inv}(w^J)$, yielding $u^J\leq_L w^J$.

The map $(u^J,u_J)\mapsto u$ is known to be a bijection between $\mathfrak{B}_n^J\times (\mathfrak{B}_n)_J$ and $\mathfrak{B}_n$, so $\varphi$ is injective.
Since $\Lambda^{L}_{w^J}\subseteq W^J$ and $\Lambda^{L}_{w_J}\subseteq W_J$, we see that $\varphi$ is length-additive. This completes the proof.
\end{proof}

\begin{prop}\label{ranksymmetric}
Let $w=w_1w_2\cdots w_n\in \mathfrak{B}_n$ be minimal non-separable. If $\{|w_{n-1}|,|w_{n}|\}=\{n-1,n\}$, then $\Lambda^{L}_{w}$ is rank-symmetric and rank-unimodal.
\end{prop}
\begin{proof} Since $w\in \mathfrak{B}_n$ is minimal non-separable,
there exists $i\in [n-1]$ such that ${\rm sts}(w_iw_n)\in\{\underline{2}1,2\underline{1}\}$ by Lemma \ref{minonsepBn}\eqref{minonsepBn:item2}.
But $\{|w_{n-1}|,|w_{n}|\}=\{n-1,n\}$, so we must have either $w_{n-1}w_{n}=\underline{n}(n-1)$ or $w_{n-1}w_{n}=n\underline{(n-1)}$. Suppose we are in the first case,
the proof for the second case is similar and will
be omitted.

Let $J=\Delta\backslash\{\alpha_{n-2},\alpha_{n-1}\}$.
Then, by Lemma \ref{varphilength-additive-bijectionmap}, we have $F(\Lambda^{L}_w)=F(\Lambda^{L}_{w^J})F(\Lambda^{L}_{w_J})$.
Since $\ell(w^J)=2n-2$, it follows from Lemma \ref{lem:uniqreducedexpre} that
$$F(\Lambda^{L}_{w^J})=1+q+q^2+\cdots+q^{2n-2},$$
which is symmetric and unimodal.
Because $w$ is a minimal non-separable element, we see that $w_J$ is separable, and hence $\Lambda^{L}_{w_J}$ is rank-symmetric and rank-unimodal by Lemma \ref{rank-symmetric}.
Therefore,  $\Lambda^{L}_{w}$ is rank-symmetric and rank-unimodal.
\end{proof}

The hypothesis that $\{|w_{n-1}|,|w_n|\}=\{n-1,n\}$ in Proposition \ref{ranksymmetric} is necessary.
Otherwise, it is difficult to predict the rank-symmetry of the set $\Lambda^{L}_{w}$.
For example, although the elements $\underline{2}3451$ and $2345\underline{1}$ are minimal non-separable in $\mathfrak{B}_5$ by Lemma \ref{minonsepBn}, the polynomial
$$\Lambda_{\underline{2}3451}^L (q)=1+2q+2q^2+2q^3+q^4+q^5$$ is not symmetric, while the polynomial  $$\Lambda_{2345\underline{1}}^L (q)=1+q+q^2+q^3+q^4+q^5$$ is symmetric.

Wei \cite{We12} showed that for any separable permutation $w$ the lower order ideal $\Lambda^{L}_{w}$ is rank-symmetric and rank-unimodal,
which was extended to all Weyl groups by Gaetz and  Gao \cite{GG20aam}.
Proposition \ref{ranksymmetric} provides a class of non-separable signed permutations whose lower order ideals are also rank-symmetric and rank-unimodal.
A natural question is to classify all (minimal) non-separable elements $w$ in a Weyl group such that the posets $\Lambda^{L}_{w}$ are rank-symmetric and rank-unimodal.

%\noindent
%{\bf Acknowledgements.}This work was partially supported by the National Natural Science Foundation of China (Grant No. 12471020, 12071377),the Natural Science Foundation of Chongqing (Grant No. CSTB2023NSCQ-MSX0706), and the Fundamental Research Funds for the Central Universities (Grant No. SWU-XDJH202305).

\end{document}